\documentclass[]{article}
\usepackage{xcolor}

\usepackage{amssymb,amsfonts,amsmath}%
\usepackage{latexsym}
\usepackage{graphicx,epsfig}

\newtheorem{thm}{Theorem}
\newtheorem{lem}[thm]{Lemma}
\newtheorem{pro}[thm]{Proposition}
\newtheorem{cor}[thm]{Corollary}
\newtheorem{df}[thm]{Definition}
\newtheorem{cnj}[thm]{Conjecture}
\newtheorem{ex}[thm]{Example}
\newtheorem{cl}[thm]{Claim}
\newtheorem{fact}[thm]{Fact}
\newtheorem{remark}[thm]{Remark}

\newcommand{\qed}{\hspace*{\fill}$\Box$}

\newcommand*{\proof}{\noindent {\bf Proof}.\,\,}

\title{\bf Constructing  sparsest $\ell$-hamiltonian saturated $k$-uniform hypergraphs for a wide range of $\ell\;\;$\thanks{This is a revised version of the published one. We fill a gap  in the proof of Lemma 10, where one case
(case 2 on page 10 in the current version) was overlooked; minor changes are scattered through Section 3.2 only.)}}

\author{Andrzej Ruci\'nski\thanks{Research supported by Narodowe Centrum Nauki, grant 2018/29/B/ST1/00426.}\\
\small Department of Discrete Mathematics\\[-0.8ex]
\small Adam Mickiewicz University\\[-0.8ex]
\small Pozna\'n, Poland\\
\small\tt rucinski@amu.edu.pl\\
\and
Andrzej {\.Z}ak \thanks{Research partially supported
by the Polish Ministry of Science and Higher Education.} \\
\small Faculty of Applied Mathematics\\[-0.8ex]
\small AGH University of Science and Technology\\[-0.8ex]
\small Krak\'ow, Poland\\
\small\tt zakandrz@agh.edu.pl
}

\begin{document}
\maketitle
\begin{abstract}
Given $k\ge3$ and $1\leq \ell< k$, an  $(\ell,k)$-cycle is one in which consecutive edges, each of size $k$, overlap in exactly $\ell$ vertices.
We study the smallest number of edges in $k$-uniform $n$-vertex hypergraphs which do not contain hamiltonian $(\ell,k)$-cycles, but once a new edge is added, such a cycle is promptly created.  It has been conjectublack that this number is of order $n^\ell$ and confirmed for $\ell\in\{1,k/2,k-1\}$, as well as for the upper range $0.8k\leq \ell\leq k-1$. Here we extend the  validity of this conjecture to the lower-middle range $(k-1)/3\le\ell< (k-1)/2$.

\end{abstract}

\section{Introduction}
A $k$-uniform hypergraph $H$ which we will be calling a $k$-graph, is a family of $k$-element subsets (edges) of a vertex set $V$.
 Given integers $1\leq \ell< k$, an \emph{$(\ell,k)$-cycle} is a $k$-graph which, for some $s$ divisible by $k-\ell$, consists of distinct vertices $v_1,\dots, v_s$ and $s/(k-\ell)$ edges
 $$\{v_1,\dots,v_k\},\;\{v_{k-\ell+1},\dots,v_{2k-\ell}\},\;\dots,\;\{v_{s-(k-\ell)+1},\dots,v_s,v_1,\dots, v_\ell\}.$$
An  \emph{$(\ell,k)$-path} is defined similarly.
Note that the number of vertices in an $(\ell,k)$-path  equals $\ell$ modulo  $k-\ell$.

A $k$-graph $H$ is \emph{$\ell$-hamiltonian saturated} (a.k.a. \emph{maximally non-$\ell$-hamiltonian}) if it is not
$\ell$-hamiltonian, but adding any new edge results in creating a hamiltonian $(\ell,k)$-cycle.

We are interested in the smallest possible number of edges, denoted by sat$(n,k,\ell)$, of an \emph{$\ell$-hamiltonian saturated} $k$-graph on $n$ vertices.
For graphs, Clark and Entringer \cite{CE} proved that sat$(n,2,1)=\lceil 3n/2\rceil$ for all $n\geq 52$.

As the problem for hypergraphs, introduced in \cite{k-sur,kk}, seems to be much harder,
we are quite satisfied with results estimating the order of magnitude of sat$(n,k,\ell)$.
Listing the results below, we silently assume that $n$ is divisible by $k-\ell$. It was observed in \cite{RZ},  Prop. 2.1, that for all $ k\geq 3$ and $1\leq \ell \leq k-1,$
\begin{align}\label{dolne_all}
{\rm sat}(n,k,\ell)=\Omega(n^{\ell})
\end{align}
  and conjectublack that this lower bound gives the correct order of magnitude.
  \begin{cnj}\label{Con} For all $k\ge3$ and $1\le\ell\le k-1$
  \begin{equation}\label{1ormany}
{\rm sat}(n,k,\ell)=\Theta(n^{\ell}).
\end{equation}
\end{cnj}
  In \cite{RZ,RZ3} we confirmed this conjecture for
 $\ell=1$, $\ell=k/2$, as well as for all $0.8k\leq \ell\leq k-1$,
(see  \cite{zak} for the
case $\ell=k-1$).
In \cite{RZ2} we proved a weaker general upper bound
 \begin{equation}\label{weakerb}
 {\rm sat}(n,k,\ell)=O\left(n^{\frac{k+\ell}{2}}\right)
 \end{equation}
and improved it
for some pairs $(k,\ell)$ in the range $\ell >k/2$.
In this paper, our main result sets another general bound on ${\rm sat}(n,k,\ell)$ which improves \eqref{weakerb} for every pair $(k,\ell)$
where $(k-2)/5<\ell< (k-1)/2$.

\begin{thm}\label{new}
Let $2\le\ell< (k-1)/2$ and $p=\max\left\{\ell,k-2\ell-1,\lceil k/2 \rceil-\ell\right\}$.
Then
$${\rm sat}(n,k,\ell)=O\left(n^{p}\right).$$
\end{thm}
Note that $p<(k+\ell)/2$ when $k-2\ell-1<(k+\ell)/2$ which is equivalent to $(k-2)/5<\ell<\lfloor k/2 \rfloor$.

The bound in Theorem \ref{new} is strong enough to confirm Conjecture \ref{Con} for a new, wide range of~$\ell$.
\begin{cor}\label{new exact}
If $(k-1)/3 \leq \ell < (k-1)/2$,
then
$${\rm sat}(n,k,\ell)=\Theta\left(n^{\ell}\right).$$
\end{cor}
 In particular, the smallest new cases of $(k,\ell)$ coveblack by Corollary \ref{new exact} include $ (6,2)$ and $(7,2)$.

 Our proof follows the general line of that in \cite{RZ3}, where the case $\ell=k/2$ was settled, but with significant alterations. First of all, we had to carefully blackefine  and recalculate  many parameters involved in the proof. An additional technical difficulty was that now we allow also odd values of $k$. However, the main obstacle, compablack with the construction in \cite{RZ3},
was due to the  gap between two consecutive disjoint edges on an  $(\ell,k)$-path, caused by  considering $\ell<k/2$.
To overcome this problem, among others, we had to prove new properties of the crucial function $\nu$ (see Section 2.1). 

\section{Construction} We will prove Theorem \ref{new} by constructing, for any large $N$ divisible by $k-\ell$, an $\ell$-hamiltonian saturated $k$-uniform hypergraph
on $N$ vertices and with $\Theta\left(N^{p}\right)$ edges. (From now on we use $N$, as $n$ is reserved for the order of a graph which plays a crucial role in the construction). In this section, we first define some parameters and then describe our construction.  We then present  a short proof of Theorem \ref{new}, the two ingblackients of which,  Lemmas \ref{pierwszy} and \ref{drugi}, will be proved in the last two sections.

\subsection{The function $\nu$}\label{fniu} In our proofs a pivotal role will be played by $(\ell,k)$-paths whose every edge draws at least $k-\ell+1$ vertices from the same fixed, relatively small set, while the remaining vertices come from a much larger set. To handle the maximum length of such paths we introduce the following function.

\begin{df}[function $\nu$]\label{defniu}\rm Given a positive integer
$x$,  let $U$ and $W$ be two disjoint sets with $|U|=x$ and $|W|=\infty$. Then
$$\nu(x)=\max_P|V(P)|,$$ where the maximum is taken over all $(\ell,k)$-paths $P$  (in the complete $k$-graph on $U\cup W$)
such that
\begin{equation}\label{constr}
U\subset V(P)\subset U\cup W\qquad\mbox{and}\qquad|e\cap U|\ge k-\ell+1\qquad\mbox{for all}\qquad e\in P.
\end{equation}
\end{df}
Note that $\nu(x)\ge x$ and $\nu(x)$ is a nondecreasing function of $x$ (just replace in $P$ one vertex of $W$ with a new
vertex of $U$). Since $\nu(x)$ is monotone, for any non-negative real number $z$ we can define

\begin{align}\label{mu}
\mu(z)=\max \left\{x: \nu(x)\leq z \right\} \quad\mbox{and}\quad\mu^*(z)=\mu(z)+1=\min\left\{x: \nu(x)> z  \right\}.
\end{align}

In the Appendix we prove several properties of function $\nu$ which will be heavily used throughout our proof.

\subsection{Parameters setting}
In this subsection we  define parameters and sets  to be used in our construction.
Set
\begin{align}
&N_0:=100k^{10} \label{n_0},
\end{align}
  let $N\geq N_0$ be an integer divisible by $k-\ell$, and
  \begin{equation}\label{wzor_n}
n:=\left\lfloor \frac{N}{11k^5} \right\rfloor.
\end{equation}
 It can be easily deduced from (\ref{wzor_n}) and  (\ref{n_0}) that
 \begin{equation}\label{Nnk}
 11k^5\le \frac Nn\le 11.5k^5\quad\mbox{and}\quad n\ge N/(11k^5)-1\ge9k^5.
 \end{equation}

Further, recall definitions in \eqref{mu} and set
\begin{align}\label{mocV}
& z:=\frac{N+4k^3}{n} -(3k-4\ell), \nonumber \\
&  x:=\mu \left( z \right)+2\lfloor k/2 \rfloor, \\
& x^*:=\mu^* \left(  z\right)+2\lfloor k/2 \rfloor +(k-2\ell)=x+(k-2\ell)+1.  \nonumber
\end{align}

\medskip

\bigskip

The following tight estimates of $N$ lie at the heart of our construction, which will become evident only at the conclusions of the proofs of the crucial Lemmas \ref{pierwszy} and \ref{drugi}. The proof is deferblack to the Appendix

\begin{pro}\label{PxiI} There exist $x_i \in \{x,x^*\}$, $i=1,\dots,n$, such that
for each $I\subset \{1,\dots,n\}$ with $|I|=n-1$,
\begin{align}
(3k-4\ell)n+\sum_{i\in I} \nu(x_i-2\lfloor k/2 \rfloor)+8k^4<N<(3k-4\ell)n+\sum_{i=1}^n \nu(x_i-2\lfloor k/2 \rfloor)-4k^3. \label{dwie_nier}
\end{align}
\end{pro}

\bigskip

Finally, we are ready to define the vertex set of the hypergraphs to be constructed.
Let
$\{A_i,B_i:\,i=1,\dots,2n\}$ be a family of $4n$ pairwise disjoint sets of sizes
\begin{align}\label{mocai}
|A_i|=\begin{cases} 2\left\lfloor k/2 \right\rfloor+\ell,\qquad\qquad i=1,\dots,n \\
2k-2\ell-3,\qquad\qquad i=n+1,\dots,2n, \end{cases}
\end{align}
and
\begin{align}\label{mocbi}
|B_i|=\begin{cases} x_i-2\left\lfloor k/2 \right\rfloor-\ell,\quad i=1,\dots,n \\
b_i\quad\qquad\qquad\qquad\qquad i=n+1,\dots,2n, \end{cases}
\end{align}
where the $x_i$'s are defined via Proposition \ref{PxiI}, while the $b_i$'s differ from each other by at most one and are chosen in such a way that
\begin{align}\label{ABN}
\sum_{i=1}^{2n}(|A_i|+|B_i|)=N.
\end{align}
The argument  that the $b_i$'s are well defined along with some bounds on them, as well as on the $x_i$'s is given in the Appendix.

\medskip

\subsection{Main construction} Let $G_1$ be a maximally non-hamiltonian
graph   with $V(G_1)=[n]=\{1,\dots,n\}$ and $\Delta(G_1)\le 5$. The existence of such a graph can be deduced for each $n\ge52$ from the results in \cite{Bondy} and \cite{CES} (see Cor. 2.6 in \cite{RZ}).
Our construction is based on the graph $G$ obtained from $G_1$ by attaching $n$  vertices $n+1,...,2n$ and $n$ edges  $\{i,n+i\}$, $i=1,\dots,n$, so that each new vertex has  degree one.

Fix $2\le\ell<{\color{black}(k-1)/2}$.
The desiblack $k$-graph $H$ will be defined on an $N$-vertex set
\begin{equation}\label{VU}
V=\bigcup_{i=1}^{2n} U_i, \quad\mbox{where}\quad U_i=A_i \cup B_i
\end{equation}
 and $A_i$, $B_i$ are given in the previous subsection (cf. \eqref{ABN}).

Before defining the edge set of $H$, we need some more terminology and notation.
For a graph $F$ and a set $S\subset V(F)$, denote by $F[S]$ the subgraph of $F$ induced by $S$. For two $k$-graphs $F_1$ and $F_2$ with $V(F_1)=V(F_2)$, we denote by $F_1\cup F_2$ the $k$-graph on the same vertex set  whose edge set is the union of the edge sets of $F_1$ and $F_2$.

For $S\subset V$, set
\begin{align*} tr(S)=\left\{i:S \cap U_{i} \neq \emptyset\right\},\quad tr_1(S)=tr(S)\cap[n],\quad\mbox{and}\quad\min(S)=\min\left\{i\in tr(S)\right\}.
\end{align*}
Note that $tr_1(S)\subset V(G_1)$. {\color{black} The set $tr(S)$ is sometimes called \emph{the trace of $S$.}}

Further, let $c(S)$ be the number of connected components of  $G^3[tr(S)]$,
where $G^3$ is the third power of~$G$, that is, the graph with the same vertex set as $G$ and with edges joining all pairs of distinct vertices which are at distance at most three in $G$.

We define the desiblack $k$-graph $H$  in terms of three other $k$-graphs, $H_1,\;H_2$, and $H_3$. Let

$$H_1^1=\left\{e\in {V \choose k}:\; \exists\{i,j\}\in G_1,\; tr_1(e)=\{i,j\},\; |A_i\cap e|\geq \left\lfloor k/2 \right\rfloor\text{ and }
|A_j\cap e|\geq \left\lfloor k/2 \right\rfloor  \right\},$$

$$H_1^2=\left\{e\in {V \choose k}:\text{ for some } i\in[n],\;tr(e)=\{i,n+i\},\;  |A_i\cap e|=\ell+1,\;
|A_{n+i}\cap e|= k-\ell-1  \right\},$$
and
$$H_1=H_1^1\cup H_1^2.$$
\begin{remark}\label{r}\rm
Note that when $k$ is odd, for an edge $e\in H_1^1$ one may actually have $tr(e)=\{i,j,r\}$, where $\{i,j\}\in G_1$, $|A_i\cap e|= \left\lfloor k/2 \right\rfloor$,
$|A_j\cap e|= \left\lfloor k/2 \right\rfloor$, and $r\in\{n+1,\dots,2n\}$, $|U_r\cap e|=1$. Note also that for an edge $e\in H_1^2$, we have $tr(e)=\{i,n+i\}\in G-G_1$. It follows that $H_1^1\cap H_1^2=\emptyset$.
\end{remark}

Further, let
$$H_2=\left\{e\in {V \choose k}:\left|e\cap U_{\min(e)}\right|\geq k-\ell+1\right\}.$$
Note that $H_1\cap H_2=\emptyset$. Indeed, if $e\in H_1$, then $|e\cap U_{\min(e)}|\le\lceil k/2\rceil<k-\ell+1$.

\begin{figure}
\begin{center}
\includegraphics[scale=0.4]{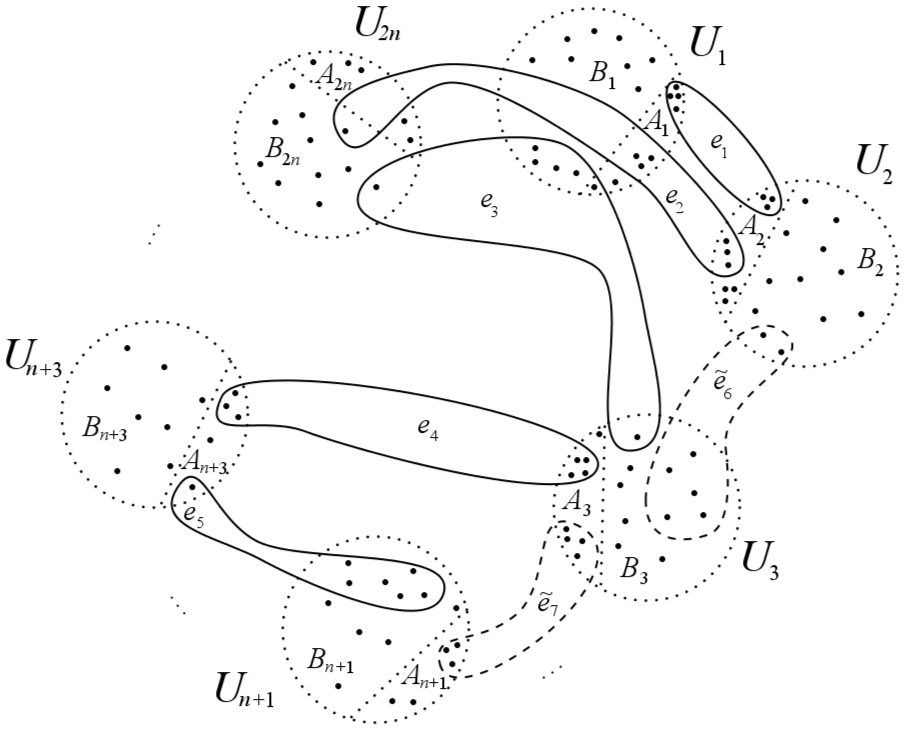}
\end{center}
\caption{Illustration of definitions of $H_1^1$, $H_1^2$, and $H_2$ for $k=7$ and $\ell=3$: {\color{black}$e_1,e_2\in H_1^1$, $e_4\in H_1^2$,  $e_3,e_5\in H_2$, while $\tilde e_6,\tilde e_7\not\in H_1\cup H_2$}.}
\label{HHH}
\end{figure}

\begin{ex}
{\rm To illustrate these definitions, let us look at Figure \ref{HHH} and the fate of the various edges depicted there. We have $k=7$ and $\ell=3$. Assume that $\{1,2\}$ is an edge of $G_1$. As $tr_1(e_1)=tr(e_1)=\{1,2\}$ and $|e_1\cap A_1|\geq |e_1\cap A_2|=3=\lfloor 7/2 \rfloor$,
 $e_1\in H_1^1$.
Further, $tr(e_2)=\{1,2,2n\}$, but $tr_1(e_2)=\{1,2\}$.
What is more,
$|e_2\cap A_1|= |e_2\cap A_2|=3=\lfloor 7/2 \rfloor$, so  $e_2\in H_1^1$ too.

Since $|e_3\cap U_1|=5=k-\ell+1$ and $\min(e_3)=1$, we have $e_3 \in H_2$. Similarly, $e_5 \in H_2$. Furthermore, $tr(e_4)=\{3,n+3\}$, $|e_4\cap A_{3}|=4=\ell+1$, and
$|e_4 \cap A_{n+3}|=3=k-\ell-1$, so
$e_4\in H_1^2$.
Finally, $|\tilde e_6 \cap U_3|=5\geq k- \ell+1$, but $\min(e_6)=2$ and $|\tilde e_6 \cap U_2|=2$.
Hence $\tilde e_6 \not\in H_1\cup H_2$. Similarly, $\tilde e_7\not\in H_1 \cup H_2$.}

\end{ex}

\medskip

Recall that
\begin{equation}\label{relacja_p_a}
p=\max\{\ell,k-2\ell-1,\lceil k/2\rceil-\ell\}.
\end{equation}
The third element of the construction is
$$H_3=\left\{e\in {V \choose k}: c(e)\leq p \right\}.$$

\begin{fact}\label{FactH3}
We have $H_1\cup H_2\subseteq H_3$.
\end{fact}
\proof If $e\in H_1$, then $|tr(e)|\le3$ and $tr(e)$ contains an edge of $G$. Thus,  $c(e)\le2\leq \ell\leq p$ and $e\in H_3$.
If $e\in H_2$, then $|e\cap U_{\min(e)}|\geq k-\ell+1$ and, consequently, $|tr(e)|\leq 1+(\ell-1)=\ell\leq p$.
Clearly, $c(e)\leq |tr(e)|$, hence $e\in H_3$ also in this case. \qed

\medskip

We are going to show (cf. Lemma \ref{pierwszy} in Section \ref{1}) that $H_1\cup H_2$ is non-$\ell$-hamiltonian.  For each $e\in\binom Vk\setminus H$, let $H+e$ be the hypergraph obtained from $H$ by adding $e$ to its edge
set.
Taking Lemma \ref{pierwszy} for granted and in view of Fact \ref{FactH3}, we define $H$ as a non-$\ell$-hamiltonian $k$-graph satisfying the containments
$$H_1\cup H_2 \subseteq H \subseteq H_3$$
and such that $H+e$ is $\ell$-hamiltonian for every $e\in H_3\setminus H$. (If $H_3$ is non-$\ell$-hamiltonian itself, we set $H=H_3$.)

\subsection{Proof of Theorem \ref{new}}

In \cite{RZ} (cf. Fact 2.2), we proved the following simple result. Let  $comp(F)$ denote the number of connected components of  a graph $F$.
\begin{cl}[\cite{RZ}]\label{bound} Let  $r$, $p$, and $\Delta$ be constants.
 If $\Delta(G)\le\Delta$, then the number of $r$-element subsets $T\subseteq V(G)$ with
 $comp(G[T])\le p$ is $O(n^p)$. \qed
\end{cl}
Theorem \ref{new} is a  consequence of Claim \ref{bound}, the construction of $H$ presented in the previous subsection, and the following two lemmas the proofs of which are deferblack to Sections \ref{1} and \ref{2}. Lemma~\ref{pierwszy}  guarantees that the definition of $H$ is meaningful.
\begin{lem}\label{pierwszy}
$H_1 \cup H_2$ is non-$\ell$-hamiltonian.
\end{lem}
On the other hand, Lemma \ref{drugi} implies quickly that $H$ is indeed $\ell$-hamiltonian saturated (see the  proof of Theorem \ref{new} below.)
\begin{lem}\label{drugi}
For every $e\in {V\choose k}\setminus H_3$, the $k$-graph $H_1\cup H_2+e$ is $\ell$-hamiltonian.
\end{lem}

\medskip

\noindent{\bf Proof of Theorem \ref{new}.} As stated in (\ref{dolne_all}), sat$(N,k,\ell)=\Omega(N^\ell)$. In order to prove the upper bound, we begin by showing that $|H|=O(N^p)$.
Observe that
$$H_3=\bigcup_{T\subset V(G)}\left\{e\in\binom V{k}:\; tr(e)=T\right\},$$
where the sum is over all subsets $T$ of $V(G)$ of size at most $k$ with $comp(G^3[T])\le p.$
Since $\Delta(G_1)\le5$, we have $\Delta(G)\le\Delta_1+1\le6$ and $\Delta(G^3)\le(\Delta_1+1)\Delta_1^2\le150$. Thus, by Claim \ref{bound} with $r\le k$, the number of such subsets $T$ is $O(n^p)$.
Moreover, by \eqref{mocV}, \eqref{ubx}, \eqref{mocai}-\eqref{mocbi} and \eqref{ubb},
\begin{equation}\label{13k5}
|U_i|=|A_i|+|B_i|\le\begin{cases}
    x_i\le x+k\le12k^5+k\le 13k^5\qquad\qquad i=1,\dots,n \\b_i+2k\le 12k^5+2k\le 13k^5\qquad\qquad i=n+1,\dots,2n\end{cases}.
\end{equation}
Hence, given $T$,
$$\Bigg|\left\{e\in\binom V{k}:\; tr(e)=T\right\}\Bigg|\le\binom{\sum_{i\in T}|U_i|}{k}\le(|T|\cdot 13k^5)^{k}=O(1).$$
 Consequently, $|H_3|=O(n^p)=O(N^p)$ and, thus, also $|H|=O(N^p)$.

It remains to show that $H$ is $\ell$-hamiltonian saturated.
 Recall that, by construction (and Lemma \ref{pierwszy}), $H$ is non-$\ell$-hamiltonian.
Let $e\in {V\choose k}\setminus H$. If $e\in H_3$ then, by the definition of $H$, $H+e$ is $\ell$-hamiltonian. On the other hand, if $e\in {V\choose k}\setminus H_3$, then $H+e\supseteq H_1\cup H_2+e$
is $\ell$-hamiltonian by Lemma~\ref{drugi}. This shows that $H$ is, indeed, $\ell$-hamiltonian saturated and  the proof of Theorem~\ref{new} is completed. \qed

\section{Proof of Lemma \ref{pierwszy}.}\label{1}

\subsection{$(\ell,k)$-paths in $H_1\cup H_2$} Before turning to the actual proof, we first establish some facts about $(\ell,k)$-paths in $H_1\cup H_2$.



\begin{fact}\label{FactH12}
If $P$ is an $(\ell,k)$-path in $H_1^2$, then $P$ has at most two edges.
\end{fact}
\proof Suppose there is an $(\ell,k)$-path $P=(e_1,e_2,e_3)$ in $H_1^2$. Then $tr(e_1)\cap tr(e_2)\neq\emptyset$ and $tr(e_2)\cap tr(e_3)\neq\emptyset$. But then, for some $j$, $tr(e_1)=tr(e_2)=tr(e_3)=\{j,n+j\}$. Since $e_1\cap e_3=\emptyset$, it follows that, in particular, $|A_{n+j}\cap e_1|=|A_{n+j}\cap e_3|=k-\ell-1$ which together exceed the size of $A_{n+j}$ set by the second part of \eqref{mocai}. \qed

\begin{fact}\label{FactH2}
If $P$ is an $(\ell,k)$-path in $H_2$, then there is an index $j\in[2n]$ such that $\min(f)=j$ for every $f\in P$, that is, every edge of $P$ draws
at least $k-\ell+1$ vertices from the same  $U_{j}$.
\end{fact}
\proof
Let $e,e'\in P$ with $|e\cap e'|=\ell$. Let $j=\min(e)$.
Since $|e\cap U_{j}|\geq k-\ell+1$, we have $|e'\cap U_{j}|\geq 1$. Hence,
$j\in tr(e')$ and so $\min(e')\leq \min(e)$. By symmetry, $\min(e)\leq \min(e')$.
Thus $\min(e')=\min(e)=j$. By transitivity, $\min(f)=j$ for every $f\in P$. \qed

\begin{cl}\label{pro2}
Let $s\geq 1$ and let $P=(e,e_1,\dots,e_s,e')$ be an $(\ell,k)$-path such that $e, e' \in H_1$ and
$e_1,\dots,e_s \in  H_2$. Then
\begin{enumerate}
\item[(i)] $\min(e_1)=\cdots =\min(e_s)\in tr_1(e)\cap tr_1(e') $;
\item[(ii)] $|\{e,e'\}\cap H_1^2|\le1$.
\end{enumerate}
\end{cl}
\proof By Fact \ref{FactH2}, $\min(e_i)=j$ for some $j\in[2n]$ and every $i=1,\dots,s$.
Since, by definition of $H_2$,  $|e_1\cap U_{j}|\geq k-\ell+1$ and  $|e_s\cap U_{j}|\geq k-\ell+1$, we have $|e\cap U_{j}|\geq 1$ and $|e'\cap U_{j}|\geq 1$ and so, $j\in tr(e)\cap tr(e')$. If, say, $e\in H_1^1$, then $tr(e)\subset[n]$, unless $k$ is odd and $|tr(e)|=3$. But then, for the unique element $r\in tr(e)\cap\{n+1,\dots,2n\}$, we have $|e\cap U_r|=1$ (cf. Remark \ref{r}), while, in fact, $|e\cap e_1|\ge2$. This means that there is $i\in tr_1(e)$ and so, $j\le i\le n$ as well.

If, on the other hand, $e,e'\in H_1^2$, then, as $tr(e)\cap tr(e')\neq\emptyset$, for some $i\in[n]$, we have $tr(e)=tr(e')=\{i,n+i\}\ni j$.
 Thus, by the definition of $H_1^2$, $|A_{n+j}\cap e|=|A_{n+j}\cap e'|=k-\ell-1$ which together exceed the size of $A_{n+j}$ set by the second part of \eqref{mocai}. This is a contradiction which excludes this case and simultaneously completes the proof of both parts, (i) and (ii).\qed

\begin{pro}\label{pro3}
Let $s\geq 1$ and $P=(e,e_1,\dots,e_s,e')$ be an $(\ell,k)$-path in $H_1\cup H_2$ such that \newline $P\cap H_1^1=\{e,e'\}$.
Then the following hold:
\begin{enumerate}
\item[(a)] $P\cap H_1^2\subset\{e_1,e_{s}\}$;

\item[(b)] If $P\cap H_1^2=\{e_1,e_{s}\}$, then $s=2$;

\item[(c)] For $i=1,\dots,s$, we have $\min(e_i)\in tr_1(e)\cap tr_1(e')$.
\end{enumerate}

\end{pro}

\bigskip

 \proof Since $s\geq 1$ and $\ell<k/2$, we have $e\cap e'=\emptyset$. If $P\cap H_1^2=\emptyset$, then the statements (a) and (b) are vacuous, while (c) follows from Claim \ref{pro2}(i).

Assume that $P\cap H_1^2 =\{f_1,\dots,f_t\} $, $t\geq 1$, where $f_i$, $i=1,\dots,t$, are listed in the order of appearance in $P$. By Claim \ref{pro2}(ii), $f_1,\dots,f_t$ are consecutive edges of $P$, while by Fact \ref{FactH12}, $t\le2$.
Recall the definition of $H_1^2$ and let $tr(f_1)=\{j,n+j\}$ for some $j\in[n]$.

When $t=2$, noticing that
$tr(f_1)\cap tr(f_{2})\neq \emptyset$ and remembering the structure of $G$, we have, in fact, $tr(f_1)=tr(f_2)=\{j,n+j\}$.
{\color{black} If $e\cap f_1 \neq \emptyset$, then $j \in tr(e)$. Indeed, otherwise $\left|e \cap U_{n+j}\right| =\left|e\cap f_1 \right| = \ell \geq 2$,
which is not possible by the definition of $H_1^1$, cf. Remark \ref{r}.  If $e \cap f_1 = \emptyset$, then,  by Claim \ref{pro2}(i) applied to the sub-path of $P$ stretching between $e$ and $f_{1}$, we have $j\in tr(e)$ too. Similar argument holds for $f_2$ and $e'$ implying that $j\in tr(e')$. Thus, $j \in tr(e)\cap tr(e')$. Since $j \leq n$, it means that $j \in tr_1(e) \cap tr_1(e')$.}

To prove (a), suppose that $e_i\in H_1^2$ for some $2\le i\le s-1$. Then, the edges $e,e_i,e'$ are pairwise disjoint. Moreover, by the definitions of $H_1^1$ and $H_1^2$,
$|A_j\cap e|\geq \left\lfloor k/2\right\rfloor$, $|A_j\cap e'|\geq \left\lfloor k/2\right\rfloor$, and $|A_j\cap e_i|=\ell+1$, which together exceed the size of $A_j$ set by the first part of (\ref{mocai}).

To prove (b), suppose that $e_1,e_s\in H_1^2$ and $s\ge 3$. Then $e_1\cap e_s=\emptyset$ and,  again by the definition of $H_1^2$,
$|A_{n+j}\cap e_1|= |A_{n+j}\cap e_s|=k-\ell-1$, which together exceed the size of $A_{n+j}$ set by the second part of (\ref{mocai}).

 It remains to prove part (c). It was already shown above that for every edge $f\in P \cap H_1^2$ we have $j=\min(f)\in tr_1(e)\cap tr_1(e')$.
 Assume now that $P\cap H_2\neq\emptyset$. Then, in view of (a) and (b), without loss of generality we may further assume that
 $e_1\in H_1^2$, while  $e_2,\dots,e_s\in H_2$.  By Claim \ref{pro2}(i) applied to the path from $e_1$ to $e'$, we conclude that for each $f\in P\cap H_2$, we have $\min(f)\in tr_1(e_1)=\{j\}$, as well
as, $\min(f)\in tr_1(e')$. Hence, $\min(f)=j\in tr_1(e)\cap tr_1(e')$ and (c) holds, indeed, for all inner edges of $P$. \qed


\subsection{Proof of Lemma \ref{pierwszy} -- the structure of phantom $C$.}

 Suppose $C$ is a hamiltonian $(\ell,k)$-cycle in $H_1\cup H_2$. We are going to show that $|V(C)|<N$ which will be a contradiction. Our proof  at some point (cf. proof of Claim \ref{cl1b}) relies on the assumption that the graph $G_1$ is not hamiltonian.

 We first consider the case when $C\cap H_1^1=\emptyset$. Then, by Fact \ref{FactH12} and Claim \ref{pro2}(ii), $C$ consists of at most two intersecting edges from $H_1^2$ and a path $P\subset H_2$. By Fact \ref{FactH2}, the bound \eqref{13k5} on $|U_j|$, and  Definition \ref{defniu} of function $\nu$ with $U=U_j$, we have, using also Proposition \ref{obs31}(b) and formula \eqref{n_0},
 $$|V(C)|\le 2k-3\ell+\nu(13k^5)\le 2k+13k^6<N_0\le N.$$

\begin{figure}
\begin{center}
\includegraphics[scale=0.35]{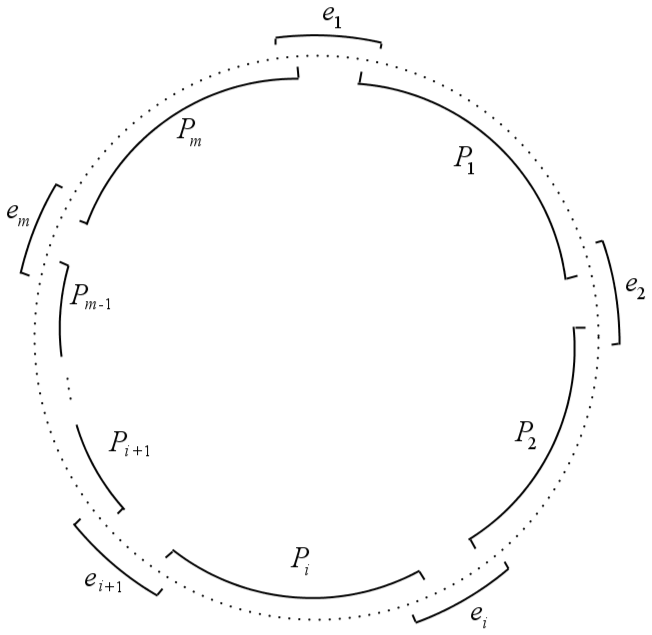}
\end{center}
\caption{{\color{black}The structure of phantom $C$.}}
\label{dawniej2}
\end{figure}

 From now on we may thus assume that $C\cap H_1^1\neq\emptyset$.
Let $M=\{e_1,\dots,e_m\}$, $m\ge1$, be a maximal set of pairwise disjoint edges of $C\cap H_1^1$, listed in the order of appearance on $C$.
Further, for  $i=1,\dots,m$, let $P_i$ be the $(\ell,k)$-path in $C$ joining the last $\ell$ vertices of $e_i$ with the first $\ell$ vertices of $e_{i+1}$, where $e_{m+1}:=e_1$. Notice that
\begin{equation}\label{C-M}
C\setminus M=\bigcup_{i=1}^m P_i,
\end{equation}
where all $P_i$'s are vertex disjoint (see Figure \ref{dawniej2}).

\begin{figure}
\begin{center}

\includegraphics[scale=0.35]{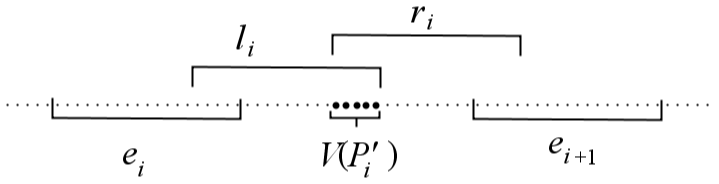}
\end{center}
\caption{{\color{black} $P'_i=\emptyset$ but $V(P'_i) \neq \emptyset$.}}
\label{nowy_rys}
\end{figure}

 Let $l_i$ be the first edge of $P_i$ and $r_i$ be the last edge of $P_i$ (note that they may coincide).
We also define $P_i'$ to be the  $(\ell,k)$-path arising from $P_i$ by removing $l_i$ and $r_i$. If $l_i=r_i$, then ${\color{black}P'_i}=\emptyset$ and $V({\color{black}P'_i})=\emptyset$. If $l_i\neq r_i$, but $P_i=\{l_i,r_i\}$, then, again, ${\color{black}P'_i}=\emptyset$, but, for the sake of the proof, we assume that $V({\color{black}P'_i})$ consists of the $\ell$ common vertices of $l_i$ and $r_i$ (see Fig. \ref{nowy_rys}).

Observe that, by the definition of $M$,
\begin{align}\label{intpi}
P_i'\subset H_1^2\cup H_2,
\end{align}
Let us now count the number $n_i$ of vertices appearing on cycle $C$ between $e_i$ and $e_{i+1}$. There are three cases.
\begin{enumerate}
\item $P_i'\neq\emptyset$: The number of vertices between $e_i$ and the beginning of $P_i'$ is exactly $k-2\ell$, and so is the number of vertices between the end of $P_{i}'$ and $e_{i+1}$. Thus, $n_i=2k-4\ell+|V(P_i')|$.
\item $P_i'=\emptyset,\; l_i\neq r_i$:   In this case, recall, $P'_i=\emptyset$ but $V(P_i')\neq\emptyset$, so the above estimates apply and, again,  $n_i=2k-4\ell+|V(P_i')|$.
\item $\l_i=r_i$: Now, $V(P_i')=\emptyset$   and  the number of vertices between $e_i$ and $e_{i+1}$ is $k-2\ell$. Thus, $n_i=k-2\ell+|V(P_i')|\le2k-4\ell+|V(P_i')|$.
\end{enumerate}

Summing up, by \eqref{C-M}, we have
\begin{equation}\label{VC}
|V(C)|{\color{black}=}mk+\sum_{i=1}^mn_i\le m(3k-4\ell)+\sum_{i=1}^m|V(P_i')|.
\end{equation}
In view of this, in order to show that $|V(C)|<N$, our plan is to utilize the left inequality in (\ref{dwie_nier}). This, in turn, will require {\color{black} us} to set strong bounds on $m$ and $|V(P_i')|$.

Beginning with the former task, recall that for each $e\in H_1^1$,
$tr_1(e)$ consists of exactly one edge of $G_1$. These edges may, however, repeat for various $e$'s, so that
 $$Tr(M):= \left\{tr_1(e): e\in M \right\} $$ is a multigraph of size $m$ on vertex set $[n]$.
 Since, for each $e \in M$ and $j\in tr_1(e)$, $|e\cap A_j|\geq \left\lfloor k/2\right\rfloor$, it follows by the first part of (\ref{mocai}) that
\begin{align}\label{delta}
\Delta(Tr(M))\leq 2,
\end{align}
 and, in particular,
\begin{align} \label{mleqn}
m\leq n.
\end{align}
To improve this bound, we distinguish between nice and problematic paths $P_i$.
Observe that each edge $e\in \left(H_1^1\cap C\right)\setminus M$ intersects some $e_i\in M$, so $e=l_i$ or  $e=r_{i-1}$.
We call an edge $l_i$ or $r_i$ \emph{bad} if it belongs to $ H_1^1$, $|P_i|\ge2$, and, resp.,  $tr_1(l_i)\neq tr_1(e_i)$
or $tr_1(r_i) \neq tr_1(e_{i+1})$.
We call $P_i$ \emph{problematic} if either $l_i$ or $r_i$ is bad, or $P_i'\cap H_1^2\neq\emptyset$. Otherwise, we call $P_i$ \emph{nice}.
{\color{black}In particular, if $P_i$ is problematic, then $|P_i|\geq 2$ and $l_i \neq r_i$}.
 Let $q$ be the number of  problematic $(\ell,k)$-paths among $P_1,\dots, P_m$.

 We next show that the presence of problematic paths makes the number of edges in $Tr(M)$ smaller.
\begin{cl}\label{cl1}
\begin{align}
m&\leq n-\frac{1}{2}\left\lceil \frac{q}{k}\right\rceil  \label{dirty}
\end{align}
\end{cl}
\proof Recall \eqref{delta}. We are going to show that problematic  paths cause some vertices to have  degrees smaller than 2 which will lead to the improvement \eqref{dirty} over \eqref{mleqn}.
Let $P:=P_i$ be problematic and assume first that there is a bad edge, say $l_i$, in $P$.
Then $tr_1(l_i)\neq tr_1(e_i)$ and, consequently, by considering separately the cases when $tr_1(l_i)\cap tr_1(e_i)=\emptyset$ and when $|tr_1(l_i)\cap tr_1(e_i)|=1$, there exists vertex $j:=j_i\in tr_1(l_i)$ such that $j\not\in tr(e_i)$ (recall Remark \ref{r}  that one might have $|tr(e_i)|=3$). Thus, by the definition of $H_1^1$, we have
$\left|\left({\color{black}l_i}\cap A_{j}\right)\setminus e_i\right|\geq \left\lfloor k/2\right\rfloor$.
Since also $|P|\geq 2$, we have $l_i\cap e_{i+1}=\emptyset$. And, obviously, by construction, $l_i$ is disjoint from all other edges in $M$. Thus, in fact,
\begin{equation}\label{=l}
\left|\left(l_i\cap A_{j}\right)\setminus (e_1\cup\cdots\cup e_{m})\right|\geq \left\lfloor k/2\right\rfloor.
\end{equation}
By symmetry, \eqref{=l} holds if $r_i$ is a bad edge of $P$.

Another reason for $P_i$ being problematic might be that $P_i'$ contains an edge $f:=f_i\in H_1^2$. Then, by the definition of $H_1^2$, there exists a vertex $j:=j_i\in tr_1(f)$ such that
$\left|f\cap A_{j}\right|= \ell+1$.
{\color{black}
Since in this case  $f$ does not intersect any edge of $M$, we have $f\cup l_i\cup r_i\subset V(P_i)$.

Thus, we may
 conclude that}, for each $i=1,\dots,m$ for which $P_i$ is problematic, there exists $j_i \in tr_1(V(P_i))$
such that
\begin{equation}\label{=l2}
\left|(V(P_i)\cap A_{j_i})\setminus (e_1\cup\cdots\cup e_m)\right|\geq
\ell+1.
\end{equation}
As $|A_{j_i}|= 2\left\lfloor k/2\right\rfloor+\ell$, inequality (\ref{=l2}) and the definition of $H^1_1$ imply that
$deg_{Tr(M)}(j_i)\leq 1$. The $j_i$'s need not be different. However, at most
$$\frac{|A_j|}{\ell+1}\le\frac{2\left\lfloor k/2\right\rfloor+\ell}{\ell+1}\le1+2\left\lfloor k/2\right\rfloor-1\le k$$
problematic paths $P_i$'s may yield the same $j$ for which $A_j$ satisfies \eqref{=l2}.
Thus, at least $\lceil q/k\rceil$ different vertices $j\in[n]$ have $deg_{Tr(M)}(j_i)\leq 1$. Therefore,
\begin{align*}
\sum_{i=1}^n deg_{Tr(M)}(i)\leq 2n-\left\lceil \frac{q}{k} \right\rceil
\end{align*}
and, consequently,
$$m=\left|Tr(M)\right|\leq n-\frac{1}{2}\left\lceil \frac{q}{k}\right\rceil. $$ \qed

\medskip

In view of Claim \ref{cl1}, we have $m\le n-1$ for $q\ge1$. Now we will get a similar improvement over $m\le n$ in the case when no problematic paths are present (unless, for some $i$, $P'_i=\emptyset$, which is, anyhow, to our advantage).

\begin{cl}\label{cl1b}
Suppose that $P_i'\neq \emptyset$ for every $i=1,\dots,m$. Then
\begin{align}
m&\leq n-1. \label{dirtyb}
\end{align}
\end{cl}
\proof If $q\geq 1$, then \eqref{dirtyb} follows by Claim \ref{cl1}. Assume that $q=0$ and suppose that $\left| Tr(M)\right|=m=n$. Then,
by (\ref{delta}), $Tr(M)$ is a 2-regular spanning subgraph of $G_1$, with possibly some parallel edge of multiplicity 2. We aim at showing that $Tr(M)$ is connected.
Since $q=0$, each $P_i$ is nice and so, by (\ref{intpi}),
$P_i' \subset  H_2$.

Let $j$ be an index guaranteed by Fact \ref{FactH2} applied to $P_i'$.
 Further, let $\bar P_i$ be the shortest extension of the path $P'_i$ within $C$ whose both end-edges belong to~$H^1_1$. Then, by Proposition \ref{pro3}(c) applied to $\bar P_i$, the traces of its end-edges contain $j\in[n]$. So, if $e_i$ is one of these end-edges, we then have $j\in tr_1(e_i)$. Otherwise, that is, when $l_i \in H_1^1$ and, thus, $l_i$ is an end-edge of $\bar P_i$,  we have $j \in tr_1(l_i)$. However, since $P_i$ is nice, $l_i$ is not
bad and so, $tr_1(e_i)=tr_1(l_i)$. Hence, $j \in tr_1(e_i)$, anyway. By symmetry,
$j \in tr_1(e_{i+1})$, too.
This means, however, that $Tr(M)$ is connected and, consequently, $Tr(M)$ is a hamiltonian cycle in $G_1$, a contradiction with the choice of $G_1$. \qed

\subsection{Proof of Lemma \ref{pierwszy} -- the length of phantom $C$.}

So far we have expressed the presumed hamiltonian $(\ell,k)$-cycle $C$ in the form \eqref{C-M} and set bounds on $m=|M|$ (see Claims \ref{cl1} and \ref{cl1b}) {\color{black} and on $|V(C)|$ (see \eqref{VC}).}
In order to take advantage of \eqref{VC},  we also need to estimate $\left|V(P_i')\right|$. We do it separately for nice and problematic paths. Recall Definition \ref{defniu} of function $\nu$ from Section \ref{fniu}.

\begin{cl}\label{cl2}
If $P_i$ is nice, then for some $j:=j_i\in[n]$,
\[\left|V(P_i')\right|\leq \nu\left(x_{j}-2\left\lfloor k/2\right\rfloor\right).\]
\end{cl}
\proof Since $P_i$ is nice, $P_i'\subset H_2$ by (\ref{intpi}). If $P_i'=\emptyset$, then the claim
{\color{black}holds by \eqref{x_nu} and \eqref{lbx}}. Let  $f\in P_i'$ and $j=\min(f)$.
Similarly, as in the proof of Claim \ref{cl1b}, we infer that  $j\in tr_1(e_i)$ and
$j \in tr_1(e_{i+1})$. 
Thus,  $|A_{j}\cap e_i|\geq \lfloor k/2 \rfloor $ and $|A_{j}\cap e_{i+1}|\geq \lfloor k/2 \rfloor$, which implies that
$\left|V(P_i')\cap U_{j}\right| \leq x_{j}-2\left\lfloor k/2\right\rfloor$.
Therefore, the claim follows by Fact \ref{FactH2} and Definition \ref{defniu} of $\nu$ with {\color{black}$U=V(P_i')\cap U_j$}.
\qed

\begin{cl}\label{cl3}
If $P_i$ is problematic, then for some $j:=j_i\in[n]$, \[\left|V(P_i')\right|\leq \nu(x_{j})+k/2.\]
\end{cl}
\proof
{\color{black} If $P_i'=\emptyset$, then the claim trivially holds. Otherwise, let $P''_i$ be the shortest extension (within $C$) of $P'_i$ with both end-edges belonging to $H_1^1$.}
By the choice of $M$, $P''_i$ exists and satisfies $P'_i \subset P''_i \subset P_i {\color{black}\cup \{e_i,e_{i+1}}\}$.
By Proposition \ref{pro3}(a,b) applied to $P''_i$, $|P''_i| \leq 4$ or $P''_i$ contains at most one edge of $H_1^2$.
In the former case the claimed inequality holds, because $|V(P'_i)|<4k$, while, by {\color{black} \eqref{x_nu} and} \eqref{lbx},
$\nu(x_j) \ge x_j \ge 10k^4$. In the latter, $P'_i$ contains
at most one edge of $ H_1^2$, as well. Moreover, this edge, if exists, is either the first or the last edge of $P_i'$. Say, it is the first. Then
the rest of $P_i'$ (i.e., $P_i'$ minus the first or the last $\ell\le k/2$ vertices) is contained in $H_2$ and either $r_i\in H_1^1$, or $r_i\in H_2$
{\color{black}(recall that since $P_i$ is problematic, $r_i \neq l_i$)}.
Hence, by Claim \ref{pro2}(i), applied to an appropriate extension of $P_i'$, there exists $j\in[n]$ such that $j=\min(f)$ for all $f\in P_i'\cap H_2$.
 Thus, $\left|V(P_i')\cap U_{j}\right| \leq |U_j|= x_{j}$ and the claim follows again by Fact \ref{FactH2} and Definition~\ref{defniu}. \qed

\bigskip

We are now in the position to finish the proof of Lemma \ref{pierwszy}.\\

\noindent
{\bf \color{black} Proof of Lemma \ref{pierwszy}.}
Let $I_1\subset [1,m]$ be the set of those indices $i$ for which 
$P'_i = \emptyset$. Let $I_2\subset [1,m]\setminus I_1$ be the set of those indices $i$ for which $P_i$ is problematic,
and $I_3=[1,m]\setminus (I_1 \cup I_2)$. Let $q_j = |I_j|$, $j \in [1,3]$.
By (\ref{VC}), Claims \ref{cl2} and \ref{cl3},  and
(\ref{iterat}) {\color{black}  applied to $x=x_{j_i}-2[k/2]$ and $t=[k/2]$,}
\begin{align*}
|V(C)|&{\color{black}\le} m(3k-4\ell)+\sum_{i\in [1,m]} |V(P'_i)|\\
&\leq m(3k-4\ell)+\sum_{i\in I_1}{\color{black}\ell}+\sum_{i\in I_2}(\nu(x_{j_i})+k/2)+\sum_{i \in I_3}\nu(x_{j_i}-2\lfloor k/2 \rfloor)\\
&\leq m(3k-4\ell)+\sum_{i\in I_1}k+\sum_{i\in I_2}(\nu(x_{j_i}-2\lfloor k/2 \rfloor)+k^2+k/2)+\sum_{i \in I_3}\nu(x_{j_i}-2\lfloor k/2 \rfloor)\\
&=m(3k-4\ell)+\sum_{i \in I_2 \cup I_3}\nu(x_{j_i}-2\lfloor k/2 \rfloor)+kq_1+(k^2+k/2)q_2 \\
&\leq m(3k-4\ell) +\sum_{i \in I_2 \cup I_3}\nu(x_{j_i}-2\lfloor k/2 \rfloor) + (k^2+3k/2)\cdot\max\{q_1,q_2\}.
\end{align*}
If $\max\{q_1,q_2\}=0$, then,  by Claim \ref{cl1b}, $m\leq n-1$  so we have $|V(C)|\le m(3k-4\ell)+\sum_{i \in [1,m]}\nu(x_{j_i}-2\lfloor k/2 \rfloor)$.
If $q_2\geq 1$, then, by Claim \ref{cl1}, $|I_2 \cup I_3|\leq m\leq n-\frac{1}{2}\left\lceil \frac{q_2}{k}\right\rceil$. 
Furthermore, $|I_2 \cup I_3| \leq n-q_1$. Hence
\begin{align*}
|I_2 \cup I_3|\leq n-\frac{1}{2}\left\lceil \frac{\max\{q_1,q_2\}}{k}\right\rceil.
\end{align*}
So, every increase of $\max\{q_1,q_2\}$ by $2k$ forces a decrease of $|I_2 \cup I_3|$ by 1. However, since by  {\color{black}  \eqref{x_nu}} and \eqref{lbx}, $\nu(x_{j_i}-2\left\lfloor k/2\right\rfloor)\geq  x_{j_i}-2\left\lfloor k/2 \right\rfloor>10k^4-k>9k^4$,  the maximum is attained when $|I_2 \cup I_3|$ is as large as possible, that is,
for $|I_2 \cup I_3|=n-1$ and $\max\{q_1,q_2\}=2k$. Hence, in either case,
\begin{align}\label{dlugoscC}
|V(C)|&\leq n(3k-4\ell)+\sum_{i \in I} \nu(x_{j_i}-2\lfloor k/2 \rfloor)+2k(k^2+3k/2),
\end{align}
where $I\subset [1,n]$ with $|I|\leq n-1$.
 Combined with the left inequality in  (\ref{dwie_nier}), this yields, with some margin, that
$|V(C)| <N$, and so $C$ cannot be a hamiltonian $(\ell,k)$-cycle, a contradiction.
\qed

\section{Proof of Lemma \ref{drugi}}\label{2}

\subsection{The idea of the proof}
In the proof of Lemma \ref{pierwszy} we supposed that there was a hamiltonian $(\ell,k)$-cycle $C$ in $H_1\cup H_2$ and got a contradiction by showing that it would be too short to cover all $N$ vertices. Now, we have at disposal just one more edge $e$ which, however, will make all the difference. In fact, despite the opposite goals these two proofs bear some similarities.

In the former proof we represented $C$ as a concatenation of several paths in $H_2$ joint together via short paths centeblack at edges of $H_1^1$.
 A crucial ingblackient of that proof was to show that there are no more than $n-1$ disjoint edges in $H_1^1\cap C$, causing the whole cycle to be too short.


 Now, we will turn that idea around and \emph{construct} a hamiltonian $(\ell,k)$-cycle in $H_1\cup H_2+e$, by
  constructing   $n$  disjoint $(\ell,k)$-paths $P_1,\dots P_n$ in $H_2$ and
 joining them by  disjoint sequences of vertices $Q_0,\dots, Q_{n-1}$ (let us call them \emph{bridges} from now on), built around edges of $H_1$.
 In fact, for technical reasons, in the forthcoming proof we will first build the bridges $Q_0,\dots, Q_{n-1}$ and only then the paths $P_1,\dots,P_n$.
The reason there were less than $n$ bridges in the proof of Lemma \ref{pierwszy} was that $G_1$ was not hamiltonian.
On the other hand, $G_1$ is \emph{maximally} non-hamiltonian and the new edge $e\not \in H$  will bring about
the missing bridge ($Q_0$). This will be done by a clever choice of two vertices  of $tr(e)$.

\bigskip

\subsection{The choice of $i$ and $j$.}\label{ij}
Let us fix $e\in {V\choose k}\setminus H_3$. Recall that, by the definition of $H_3$, $c(e)\geq p+1$, where $p$ was defined in \eqref{relacja_p_a}.
We are going to choose carefully two vertices, $i$ and $j$, in $tr(e)$. They have to come from different components of $G^3[tr(e)]$. In particular,  $ij\not\in G$. Even more,
if $i=n+i'$ or $j=n+j'$ for some $1\le i',j'\le n$, then also, respectively, $ij', i'j,i'j'\not\in G_1$. (This is, in fact, why we consideblack components in $G^3[tr(e)]$, and not just in $G[tr(e)]$.) The bottom line is that, due to being maximally non-hamiltonian, $G_1$ possesses a hamiltonian path connecting $i$ (or its unique neighbor) with $j$ (or its unique neighbor). We will ultimately build a  hamiltonian $(\ell,k)$-cycle in $H_1\cup H_2+e$ by following this path in $G_1$.

Let $C_1,\dots,C_r$
be connected components of $G^3[tr(e)]$. Further, let
$$\rho(C_t)=\max\{|e\cap U_v|:v\in V(C_t)\},\quad t=1,\dots,r.$$ Without loss
of generality we may assume that
$$\rho(C_1)\geq \rho(C_2)\geq \cdots \geq \rho(C_r).$$

We now choose $i$ and $j$.
If $\rho(C_1)\leq \ell$, then  $i=\min(e)$. Otherwise, let $i\in V(C_1)$ be such that
$$|e\cap U_i|=\rho(C_1)\geq \ell+1.$$  Let $X$ be the vertex set of this component of $G^3[tr(e)]$ which contains vertex $i$ (e.g., $X=V(C_1)$ in the latter case)
and let
$Y=tr(e)\setminus X$. Set
$$e_X=e\cap\bigcup_{v\in X}U_v\quad\mbox{and}\quad e_Y=e\cap\bigcup_{v\in Y}U_v.$$ Clearly,
\begin{align}\label{exey}
e=e_X\cup e_Y.
\end{align}
Further, if $\rho(C_2)\leq \ell$, then $j=\min(e_Y)$. Otherwise, let
$j\in V(C_2)$ be such that
$$|e\cap U_j|=\rho(C_2)\geq \ell+1.$$ Note that in the latter case $X=V(C_1)$, so, indeed, $i$ and $j$ always belong to different components of $G^3[tr(e)]$.

\bigskip

Now we establish upper bounds on the cardinalities of some parts of $e$.
Since $c(e)\geq p+1$,
\begin{equation}
|e\cap U_t |\leq k-p \text{ for every } t\in tr(e), \label{cel1}
\end{equation}
and, in particular,
\begin{equation}
|e_X|\leq k-p.\label{cel2}
\end{equation}
Note that, by (\ref{exey}) and (\ref{cel2}), we also have $e(Y)\ge p$.
Inequality \eqref{cel1} can be improved in most cases.
\begin{fact}\label{ute}
If $t\in tr(e)\setminus \{i,j\}$, then
\begin{align*}
|e\cap U_t|\leq \ell.
\end{align*}
\end{fact}
Proof. If $\rho(C_1)\leq  \ell$ then the claim is obvious.
Suppose $\rho(C_1)\geq  \ell+1$. Thus, $|e\cap U_i |\geq \ell+1$.
If $t\in X\setminus\{i\}$, then, by (\ref{cel2}) and  (\ref{relacja_p_a}),
\begin{align*}
|e\cap U_t |\leq k-p-|e\cap U_i|\leq k-p-(\ell+1)\leq \ell.
\end{align*}
 Let $t\in tr(e)\setminus X=Y$.
If $\rho(C_2)\leq  \ell$, then, again, the claim is obvious.
So, suppose $\rho(C_2)\geq  \ell+1$. Hence, $|e\cap U_j|\geq \ell+1$.
Note that since $|tr(e)|\geq c(e)\geq p+1$, we have $|tr(e)\setminus\{i,j,t\}|\geq p-2$, and so
\begin{align*}
|e\cap (U_i \cup U_j \cup U_t)|\leq k-p+2.
\end{align*}
Thus, again by (\ref{relacja_p_a}),
\begin{equation*}
|e\cap U_t |\leq k-p+2-|e\cap U_i|-|e\cap U_j |\leq k-p+2-2(\ell+1)\leq 1<\ell.
\end{equation*} \qed

\noindent
\subsection{Construction of bridge $Q_0$.}
The construction of $Q_0$ is based on the extra edge $e$ and the choice of $i$ and $j$ from $tr(e)$.  Let us order the vertices of $e$ so that, going from left to right, it begins with all vertices of $e\cap U_j$, followed by all remaining vertices of $e(Y)$. Symmetrically, going from right to left, it begins with all vertices of $e\cap U_i$, followed by  the remaining vertices of $e(X)$.

We first we construct an $(\ell,k)$-path $Q'_0$ which
is the main part of $Q_0$.
We consider four cases with respect to $i$ and $j$, which, owing to symmetry, blackuce to just two (with two further subcases in one of them).

\smallskip

\noindent
{\bf Notation for diagrams.} The forthcoming constructions will be illustrated by diagrams in which  the following
notation is applied. Recall that for each $s=1,\dots,2n$, $U_s=A_s\cup B_s$. Any vertex of $A_s$ will be represented by the symbol $a_s$. Similarly, $b_s$ will stand for any vertex of $B_s$, while $u_s$ for any vertex of $U_s$. The asterisk $*$ will fill in for any vertex of $V=\bigcup_{s=1}^{2n}U_s$, or, on one occasion, of $\bigcup_{s=n+1}^{2n}B_s$. Moreover, all vertices appearing in the diagrams will be distinct.

\bigskip

Suppose first that $i,j\in \{1,\dots, n\}$.
Let $Q'_0$ be a 3-edge $(\ell,k)$-path with the edge $e$ in the middle and two edges $e'$ and $e''$ from $H_2$.
The first edge $e'$  of $Q'_0$ begins with $k-\ell$ vertices of $B_j$ and ends with the first $\ell$ vertices of $e$, while the last (third) edge $e''$ of $Q'_0$ begins with the last $\ell$ vertices of $e$ and ends with $k-\ell$ vertices of $B_i$ (see  diagram \eqref{eee} below).

\begin{equation}\label{eee}
Q'_0=\underbrace{b_j\dots b_j}_{k-\ell} \underbrace{\overbrace{u_j**}^{e_Y}\overbrace{*u_i}^{e_X}}_{e}\underbrace{b_i\dots b_i}_{k-\ell}.
\end{equation}
Recall that either $j=\min(e_Y)$ or $|U_j\cap e|\geq \ell+1$.
Consequently, in each case $\min(e')=j$ and $|e'\cap U_j|\ge k-\ell+1$, so  $e'\in H_2$. Similarly, $e''\in H_2$. 

If $i=n+i'$, then we modify the right end of $Q'_0$ as follows.  If $|e\cap A_i|\le k-\ell-2$, then we replace the last $\ell$ vertices of $e''$ with $k-\ell-1$ vertices of $A_i$, followed by $\ell+1$  vertices of $A_{i'}$ (see the R-H-S of diagram \eqref{eeeee}).

 \begin{equation}\label{eeeee}
Q'_0=\underbrace{b_j\dots b_j}_{k-\ell}\underbrace{\overbrace{u_j**}^{e_Y}\overbrace{*u_i}^{e_X}}_{e}
\underbrace{b_i\dots b_i}_{k-2\ell}\underbrace{a_i\dots a_i}_{k-\ell-1}\underbrace{a_{i'}\dots a_{i'}}_{\ell+1}.
 \end{equation}
 This way, edge $e''$ is replaced by  edges $e''_1\in  H_2$ and  $e''_2\in H_1^2$. Since $|e\cap A_i|\le k-\ell-2$, we have, indeed, at least $(2k-2\ell-3)-(k-\ell-2)=k-\ell-1$ vertices of $A_i$ available. (As for $A_{i'}$, by \eqref{mocai}, $|A_{i'}|\ge k-1+\ell$, and only at most $k-2$ vertices of $A_{i'}$ may belong to $e$.)

If $|e\cap A_i|\geq k-\ell-1$, we modify $Q'_0$ as indicated in the R-H-S of diagram \eqref{eeee}.
 \begin{equation}\label{eeee}
Q'_0=\underbrace{b_j\dots b_j}_{k-\ell}
\underbrace{\overbrace{u_j**}^{e_Y}\overbrace{*a_i\dots a_i}^{e_X}}_{e}
\underbrace{a_i\dots a_i}_{k-2\ell-1}\underbrace{a_{i'}\dots a_{i'}}_{\ell+1}.
 \end{equation}
Note that now, again, we have just one edge to the right of $e$ and this is an edge of $H_1^2$. Furthermore, by (\ref{cel1}) and (\ref{relacja_p_a}),
\begin{align*}
| Q'_0\cap A_i|\leq k-p+k-2\ell-1\leq 2k-2\ell-3,
\end{align*}
so, this construction is feasible.

The case $j=n+j'$  is analogous. In summary, depending on the case, the path $Q_0'$ consists of three to five edges, all contained in $H_1\cup H_2+e$.
{\color{black}
To simplify further notation, from now on, let us assume (w.l.o.g.) that $i\in\{1,n+1\}$ and $j\in\{n,2n\}$.
In fact, we may  arbitrarily renumber vertices $1,\dots,n$ and, accordingly, vertices $n+1,\dots, 2n$.
Since in the rest of the construction
we are  going to  use only edges $e'$ of $H_2$ that intersect exactly one of the sets $U_i$ with $1\leq i \leq n$, such a renumbering will not affect
the sets $U_{\min(e')}$ (which are crucial for the edges of $H_2$), regardless of how may sets $U_i$ with $n+1\le i\le 2n$ are intersected by $e'$.}

We complete the construction of $Q_0$ by adding $k-2\ell$ new vertices from $B_{n}$ on the left of $Q'_0$
and $k-2\ell$ new vertices from $B_{1}$ on the right of $Q'_0$, that is,
\begin{align}\label{diagram_q0}
Q_0 = \underbrace{b_n \dots b_n}_{k-2\ell}Q'_0\underbrace{b_1 \dots b_1}_{k-2\ell}.
\end{align}
Note that $k-2\ell\ge1$ and that   $Q_0$ always begins with at least $k-\ell+1$ vertices from $U_{n}$ and ends with at least $k-\ell+1$ vertices from $U_{1}$. Also, technically, $Q_0$ is not an $(\ell,k)$-path as at either end it is, on purpose, ``unfinished''.

\medskip

 Before continuing with the construction, let
us  summarize how many vertices have been taken by $Q_0$ from each set $A_t$, $t\in[n]$.
To this end, let us partition the set $[n]$ into two subsets
\begin{align}\label{T1T2}
T_1&=\{t\in[n]:t\not\in tr(e) \text{ and } n+t \not\in tr(e) \}, \nonumber \\
T_2&=[n]\setminus T_1
\end{align}
and observe that
\begin{align}\label{T2}
T_1\subseteq[2,n-1]\quad\mbox{and}\quad|T_2|\leq |tr(e)|\le k.
\end{align}
Trivially, by the construction of $Q_0$, for all $t\in T_1$,
\begin{align}\label{T1}
(U_t \cup U_{n+t})\cap Q_0=\emptyset.
\end{align}

\begin{fact}\label{pezero}
\begin{align}\label{lp0}
| Q_0\cap A_t|\leq \begin{cases} k-p \qquad\mbox{ for }\quad t\in  \{1,n\}\;, \\
\ell \qquad\qquad\mbox{ for }\quad t\in T_2\cap[2,n-1]\;,\\
0 \qquad\qquad\mbox{ for }\quad t\in T_1 \;.
\end{cases}
\end{align}
\end{fact}
Proof. If $t\in T_1$ then the statement follows from (\ref{T1}).
If $t\in T_2\cap [2,n-1]$, then by the construction of $Q_0$,
\begin{align*}
Q_0\cap A_t\subseteq Q_0\cap U_t=e\cap U_t
\end{align*}
and the second line of (\ref{lp0}) holds by Fact \ref{ute}.

Let $t=1$. If $i=1$, then the R-H-S of $Q'_0$ is like
in diagram (\ref{eee}), and so, by (\ref{cel1}),
$$| Q_0\cap A_1|=|e\cap A_1|\leq |e\cap U_1 |\leq k-p.$$  If, on the other hand,
$i=n+1$, then consider two cases with respect to whether $1\in tr(e)$ or not.

If $1\not\in tr(e)$, then by diagrams (\ref{eeeee}) or (\ref{eeee}),   and by (\ref{relacja_p_a}),
\begin{align*}
| Q_0\cap A_1|=\ell+1\leq k-p.
\end{align*}
(To see the last inequality one has to check all 3 cases for $p$.)

On the other hand, if $1\in tr(e)$, the procedure selecting $i$ implies that
$$|e\cap U_{n+1}|=\rho(C_1)\ge\ell+1.$$
Furthermore, as $1$ and $n+1$  are two vertices of the same component of $G$, and thus of $G^3$, we have $\{1,n+1\}\subseteq X$ and, by \eqref{cel2},
\begin{equation}\label{ek-p}
|e\cap U_1|+|e\cap U_{n+1}|\le |e_X|\le k-p.
\end{equation}
Hence, again by diagrams (\ref{eeeee}) or (\ref{eeee}),
$$|Q_0\cap A_1|\le|e\cap U_1|+(\ell+1)\le|e\cap U_1|+|e\cap U_{n+1}|\le k-p.$$

 The proof for $t=n$ is analogous, except that in the case $j=2n$, $n\in tr(e)$, to get an analog of (\ref{ek-p}), instead of \eqref{cel2} we use the inequality $|tr(e)\setminus\{n,2n\}|\ge c(e)-1\ge p$ which immediately implies that
  $$|e\cap U_n|+|e\cap U_{2n}|\le k-p.$$ \qed

\subsection{Construction of bridges $Q_1,\dots,Q_{n-1}$.}
Since $G_1$ is maximally non-hamiltonian and $1n\not\in G_1$, there is a hamiltonian path in $G_1$ which begins at vertex $1$ and ends at vertex $n$.
W.l.o.g., we assume that its vertex sequence is $1,2,3,\dots,n-1,n$.
Based on this hamiltonian path we will  build a hamiltonian $(\ell,k)$-cycle in $H$.

First, we construct $n-1$ pairwise disjoint edges, $e_1\dots,e_{n-1}\in H_1$, such that they are also disjoint from $e$ and for each $t=1,\dots,n-1$,  $e_t$ contains $\lfloor k/2 \rfloor$ vertices from $A_{t}$ followed, if $k$ is odd, by one  vertex from $\bigcup_{s=n+1}^{2n}B_s$  and
then $\lfloor k/2 \rfloor$ vertices from $A_{t+1}$ (see the diagram below).
\[e_t=\underbrace{a_t\dots a_t}_{\lfloor k/2 \rfloor}(*)\underbrace{a_{t+1}\dots a_{t+1}}_{\lfloor k/2 \rfloor}.\]
Thus, for each $s=2,\dots, n-1$ we need $2\lfloor k/2 \rfloor$ vertices of $A_s$
which is feasible by (\ref{lp0}) and (\ref{mocai}), while for $s\in\{1,n\}$ we only need $\lfloor k/2 \rfloor$ vertices of $A_s$, which is again possible by (\ref{lp0}) and (\ref{mocai}), and the definition of $p$ in (\ref{relacja_p_a}).

Next we set aside pairwise disjoint $(k-2\ell)$-element sequences of vertices $L_1,\dots,L_{n-1}$ and $R_1,\dots,R_{n-1}$ which are also disjoint from $Q_0\cup e_1\cup\cdots\cup e_{n-1}$ and such that for all $ t=1,\dots n-1$ we have $L_t \subset B_{t}$, while
\begin{align}
&R_t \subset A_{n+{t+1}}\quad \text{ if }\quad t+1\in T_1\;, \nonumber\\
&R_t \subset B_{t+1} \qquad\text{ if }\quad t+1\in T_2\;, \nonumber
\end{align}
which is feasible by (\ref{mocai}) together with (\ref{T1}), and (\ref{mocbi}) together with (\ref{lbx}), and the bound $|Q_0\cap B_t|\le|Q_0|<7k$.
Finally, for all $ t=1,\dots n-1$ set
\begin{align*}
Q_t =L_t,e_t,R_t,
\end{align*}
 that is,
\begin{equation}\label{Qt}
Q_t =\begin{cases}  \underbrace{b_t\dots b_t}_{k-2\ell}\underbrace{a_t\dots a_t}_{\lfloor k/2 \rfloor}(*)\underbrace{a_{t+1}\dots a_{t+1}}_{\lfloor k/2 \rfloor}
 \underbrace{a_{n+t+1}\dots a_{n+t+1}}_{k-2\ell} \qquad\qquad\text{ if }\quad t+1\in T_1\;,\\
\underbrace{b_t\dots b_t}_{k-2\ell}\underbrace{a_t\dots a_t}_{\lfloor k/2 \rfloor}(*)\underbrace{a_{t+1}\dots a_{t+1}}_{\lfloor k/2 \rfloor}
 \underbrace{b_{t+1}\dots b_{t+1}}_{k-2\ell}\; \qquad\qquad\qquad\text{ if }\quad t+1\in T_2\;.
\end{cases}
\end{equation}
\\

\noindent
So far we have constructed all bridges. Let us summarize how many vertices of each set $U_t$, $t\in[n]$, were consumed by them. In addition, for future purposes, we are also interested in the usage of $A_{n+t}$,  $t\in T_1$.
Let $Q=\bigcup_{t=0}^{n-1} Q_t$ (here $Q_t$'s are understood as sets, not sequences).
\begin{fact}\label{qi} We have the following bounds.
\begin{enumerate}
\item[(i)] For each $t\in T_1$\;,  $|Q\cap A_t|=2\lfloor k/2 \rfloor$, $|Q\cap B_t|=k-2\ell$, and  $|Q\cap A_{n+t}|= k-2\ell$.
\item[(ii)] For each $t\in T_2$\;,  $|Q\cap U_t|\le2\lfloor k/2 \rfloor+4k$.
\end{enumerate}
\end{fact}
\proof In general, $Q\cap U_t=(Q_0\cap U_t)\cup (Q_t\cap U_t)\cup (Q_{t-1}\cap U_t)$, where we assume $Q_n=\emptyset$ for convenience. By (\ref{T1}), when $t\in T_1$, we have $Q_0\cap U_t=\emptyset$ and
$Q_0 \cap A_{n+t}=\emptyset$. Also then, by inspecting \eqref{Qt}, $|Q_t\cap A_t|=\lfloor k/2 \rfloor$ and $|Q_t\cap B_t|=k-2\ell$, while $|Q_{t-1}\cap A_t|=\lfloor k/2 \rfloor$, $|Q_{t-1}\cap B_t|=0$ and $|Q_{t-1}\cap A_{n+t}|=k-2\ell$. This proves part (i).

When $t\in[2,n-1]\cap T_2$, we have $|Q_0\cap U_t|=|e\cap U_t|\le \ell$ by Fact \ref{ute}, and, again by inspection, $|Q_t\cap U_t|=|Q_{t-1}\cap U_t|=\lfloor k/2 \rfloor+k-2\ell$, so, altogether, $|Q\cap U_t|\le2\lfloor k/2 \rfloor+2(k-2\ell)+\ell\le 2\lfloor k/2 \rfloor+4k$.

Consider now the case $t=1$. Then $i=1$ or $i=n+1$.
If $i=1$, then bounding trivially $|e\cap U_1|\le k$, by  \eqref{eee} and \eqref{diagram_q0}, we have $|Q_0\cap U_1|\le k+(k-\ell)+(k-2\ell)$. This, together with $|Q_1\cap U_1|=\lfloor k/2 \rfloor+k-2\ell$, yields that
$$|Q\cap U_1|\le \lfloor k/2 \rfloor +4k-5\ell\le \lfloor k/2 \rfloor +4k.$$
If $i=n+1$, then by \eqref{eeeee}, \eqref{eeee} (with $i'=1$), \eqref{diagram_q0} and \eqref{Qt}, and again bounding  $|e\cap U_1|\le k$, we obtain
$$|Q\cap U_1|\le (\ell+1)+(k-2\ell) + (\lfloor k/2 \rfloor + k - 2\ell) +4  = \lfloor k/2 \rfloor +3k -3\ell+1\le \lfloor k/2 \rfloor +4k.$$
The case $t=n$ is very similar.
\qed

\medskip

\noindent
\subsection{Construction of paths $P_1,\dots,P_n$.}\label{conP}
Next, we construct $n$ {\color{black}pairwise} vertex disjoint $(\ell,k)$-paths $P_t\subseteq H_2$, $t=1,\dots,n$, such that
each $P_t$ consists of all vertices from $U_{t}\setminus Q$
and some vertices from $\bigcup_{s=n+1}^{2n} U_s\setminus Q$, so that together with the sequences $Q_0,\dots,Q_{n-1}$ they exhaust all $N$ vertices and, after some mending, will yield the ultimate hamiltonian $(\ell,k)$-cycle.

By the definition of $H_2$ and Fact \ref{FactH2}, each edge $f\in P_t$ will have to satisfy $\min(f)=t$ and   $|f\cap (U_{t}\setminus Q)|\ge k-\ell+1$.
We are going to build the paths $P_1,\dots,P_t$, in two stages.

\subsubsection*{Abstract Construction}
First, instead of $\bigcup_{s=n+1}^{2n} U_s$, we use vertices from some (abstract and disjoint from $V$) infinite set $W$
and  construct paths $P'_1,\dots,P'_t$ which are as large as possible and each edge $f\in P'_t$  satisfies    $|f\cap (U_{t}\setminus Q)|\ge k-\ell+1$. By Definition \ref{defniu} of function $\nu$ with $U=U_t\setminus Q$ we have $|V(P_t')|=\nu(|U_t\setminus Q|)$.
It will turn out that the total length of these paths and the sequences $Q_0,\dots,Q_{n-1}$ exceeds $N$, so in the second stage
 we will truncate them to the total length $N$ (by removing some vertices of $W$) and, finally, replace the remaining vertices of $W$ by those in $\bigcup_{s=n+1}^{2n} U_s$, obtaining the desiblack paths $P_1,\dots,P_t$.

We first estimate the lengths of the paths $P'_1,\dots,P'_t$. By Fact \ref{qi}(i), \eqref{mocai}, and \eqref{mocbi}, for $t\in T_1$ we have $|U_t\setminus Q|=x_{t}-\left(2\lfloor k/2 \rfloor+k-2\ell\right)$. Thus,  by (\ref{mocV_wn}) and (\ref{mocV_wn1}),
\begin{equation}\label{dlpiT1}
|V(P'_t)|=\nu\left((x_{t}-2\lfloor k/2 \rfloor)-(k-2\ell)\right)=\nu(x_{t}-2\lfloor k/2 \rfloor) \text{ if } t \in T_1
\end{equation}

Similarly (but understandably with less precision), by Fact \ref{qi}(ii), \eqref{mocai},  \eqref{mocbi}, and  (\ref{iterat}),  we have
\begin{equation}\label{dlpiT2}
|V(P'_t)| {\color{black}\geq} \nu\left((x_{t}-2\lfloor k/2 \rfloor) -4k\right)\geq\nu(x_{t}-2\lfloor k/2 \rfloor) -4k^2\text{ if } t\in T_2.
\end{equation}

Notice that $|Q_t|=3k-4\ell$ for all $t=1,\dots,n-1$ and, as $Q_0'$ has at least 3 edges, $|Q_0|\ge 2(k-2\ell)+3(k-\ell)+\ell\ge 3k-4\ell$. Using these estimates and recalling (\ref{T1T2}), (\ref{T2}), \eqref{dlpiT1}, and (\ref{dlpiT2}), we now bound from below the total number $N'$ of vertices appearing in all so far constructed objects.

\begin{align*}
&N'=\sum_{t=0}^{n-1}|Q_t|+\sum_{t=1}^n |V(P_t')|\\
&\geq (3k-4\ell)n+\sum_{t\in T_1} \nu(x_t-2\lfloor k/2 \rfloor) +\sum_{t\in T_2}(\nu(x_t-2\lfloor k/2 \rfloor)-4k^2)\\
&\geq (3k-4\ell)n+\sum_{t=1}^n\nu(x_t-2\lfloor k/2 \rfloor)-4k^3>N,
\end{align*}
where the last inequality holds by  (\ref{dwie_nier}).

\subsubsection*{Trimming}
\begin{figure}
\begin{center}
\includegraphics[scale=0.4]{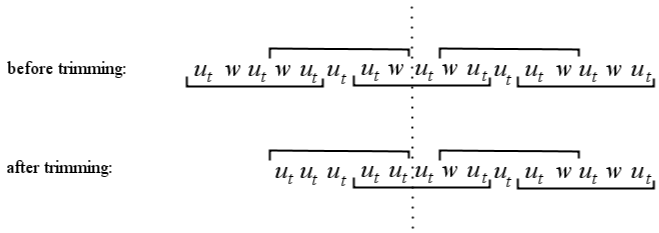}
\end{center}
\caption{\color{black} Illustration of trimming for $k=5$ and $\ell=2$;
the segment to the right of the dotted line remains unchanged, while the one to the left retains only of vertices from $U_t$.}
\label{trim_fig}

\end{figure}

 Recall that $N$ is divisible by $k-\ell$. It is easy to check that the same is true for $N'$.
 As long as $N'>N$ we apply the following iterative procedure of trimming the paths $P'_1,\dots,P'_t$: choose a path, which currently contains
{\color{black}the largest number} of vertices of $W$ and remove from it \emph{precisely}
 $k-\ell$ leftmost vertices  of  $W$ (according to the order of their appearance on the path).
 As, by \eqref{mocai} -- \eqref{ABN}, \eqref{lbb}, \eqref{VU} and Fact \ref{qi}
\begin{align}\label{wcappttotal}
\Big|\bigcup_{t=1}^n \left( V(P_t') \cap W \right)\Big| &\geq N'-{\color{black}\sum_{t=0}^{n-1}|Q_t|}-\sum_{t=1}^n|U_t|>N-{\color{black}5kn}-\sum_{t=1}^n|U_t| \nonumber \\
&=\sum_{t=n+1}^{2n}|U_t|-{\color{black}5kn}\ge n\cdot \min b_t-{\color{black}5kn} \ge (4k^4{\color{black}-5k})n,
\end{align}
a  path with at least $k-\ell$ vertices of $W$ exists (as long as $N'>N$).
It is easy to see that, treating the remaining vertices of the truncated path as consecutive, we  obtain a new, shorter (by $k-\ell$) path such that each of its edges still has at least $k-\ell+1$ vertices of $U_{t}\setminus Q$, {\color{black}see Fig. \ref{trim_fig}.
Indeed, the edges to the right of the rightmost removed element (doted line in Fig. \ref{trim_fig}) remain the same as before trimming (due to the fact that we have removed exactly $(k-\ell)$ leftmost vertices of $W$), while those to the left  have now all vertices in $U_{t}\setminus Q$. For the remaining edge (the one with vertices to the left and to the right) we argue similarly. Its part to the right remains
unchanged (and so has the same number of vertices from $U_t \setminus Q$ as before trimming), while the part to the left has now all vertices in $U_{t}\setminus Q$ (at least as many as before trimming).}

 We conclude the procedure when the current number of vertices in all the paths and sequences $Q_0,\dots,Q_{n-1}$ (which remain untouched) reaches $N$. Let the resulting paths be denoted by $P''_1,\dots,P''_n$.

{\color{black}
Furthermore,
note that by \eqref{dlpiT1}, \eqref{dlpiT2}, \eqref{x_nu} and \eqref{lbx}, at the beginning of the trimming we had
\begin{align}\label{wcappt}
\left| V(P'_t) \cap W \right| &= \left| V(P'_t)\setminus U_t\right| \geq \nu(x_t-k)-4k^2-x_t \geq \frac{k+1}{k}(x_t-k)-4k^2-x_t  \nonumber \\
&=\frac{x_t}{k}-4k^2-k \geq 10k^3-4k^2-k \geq 5k^3.
\end{align}
Since at every stage we removed vertices from a path with the largest number of vertices in $W$, by \eqref{wcappttotal} and \eqref{wcappt},
\begin{align}\label{wcapptfinal}
\left| V(P''_t) \cap W \right| \geq \min\{5k^3, 4k^4-5k -(k-\ell)\} = 5k^3.
\end{align}}

  \subsubsection*{Eradicating}
  We still have to eradicate the remaining vertices of $W$, that is, to replace them by the vertices of $\bigcup_{s=n+1}^{2n} U_s$. While doing so, we will also prepare the structure of the paths for the final concatenation into a hamiltonian $(\ell,k)$-cycle. In fact, this preparation will mostly affect only the first edge, call it $f_t''$, of $P_t''$ for $t\in T_1$.

 \textbf{Preparation:}
{\color{black}We first change the order of the first $k$ vertices of $P''_t$, so that the vertices on positions $\ell +1,\ell+2,\dots k$ are all from $U_t$. This is possible because $f''_t$ (as well as every other edge of $P''_t$)
contains at least $k-\ell+1$ vertices from $U_t$.}
Note that this operation may also affect the second edge of $P_t''$, but it will still have at least $k-\ell+1$ vertices from $U_t$.
The remaining edges of $P''_t$, as disjoint from $f''_t$, remain unchanged.
Let us call the resulting path $P_t'''$ and its first edge $f_t'''$.
Focusing on $f_t'''$, we see that among its first $\ell$ vertices at least one is from $U_t$ (because $f'''_t$ has at least $k-\ell+1$
vertices from $U_t$). Now, if
there are more than one vertices like this,  we swap all but one of them  with  arbitrary vertices of {\color{black}$W\cap \left( P_t'''\setminus f'''_t \right)$
(note that by \eqref{wcapptfinal} there are enough vertices of $W$ in $P'''_t$ to do this)}.
After this operation the number of vertices from $U_t$ in every edge (but $f_t'''$) can only increase, so
still each edge  has  at least $k-\ell+1$ vertices from $U_t$.

Finally, if necessary, we move the unique vertex of $U_t$ among the first $\ell$ vertices to the $\ell$-th position and,
if it belongs to $B_t$, we exchange it  with a vertex of $A_t$ (which also belongs to $P_t'''$). Such a vertex exists, since, by Fact \ref{qi}(i), out of all vertices of $A_t$, precisely $2\lfloor k/2 \rfloor$ were used by $Q$, while the remaining $\ell$ are sitting somewhere on the path $P_t'''$. In summary, after these changes we obtain a new path $P_t''''$ such that, for each $t\in T_1$, the structure of its first edge is
\begin{align}\label{lpt}
f_t''''=\underbrace{w,\dots,w}_{\ell-1},a_t, \underbrace{u_t\dots,u_t}_{k-\ell}.
\end{align}

\textbf{Replacement:} Finally, to obtain the desiblack paths $P_t\in H_2$, we replace  the  vertices of $W$ in $\bigcup_{t=1}^nV(P_t'''')$ by the vertices of $\bigcup_{s=n+1}^{2n} U_s$ in the following order.
First, for each $t\in T_1$,  we replace the $\ell-1$ vertices of $W$ at the left end of $f_t''''$ by vertices from $A_{n+t}$. This is possible, since
by (\ref{mocai}) and Fact \ref{qi}, there are at least $k-3\ge\ell-1$ vertices of $A_{n+t}$ unused so far. As a result, the first edge of each path $P_t$, $t\in T_1$, by (\ref{lpt}), takes the form

\begin{align}\label{lpt1}
f_t=\underbrace{a_{n+t},\dots,a_{n+t}}_{\ell-1},a_t, \underbrace{u_t\dots,u_t}_{k-\ell}.
\end{align}
The remaining
vertices of $W$ in $\bigcup_{t=1}^nV(P_t'''')$ are replaced arbitrarily.
\bigskip

\noindent
\subsection{Construction of the hamiltonian cycle $C$.}\label{conC}
\begin{figure}
\begin{center}
\includegraphics[scale=0.5]{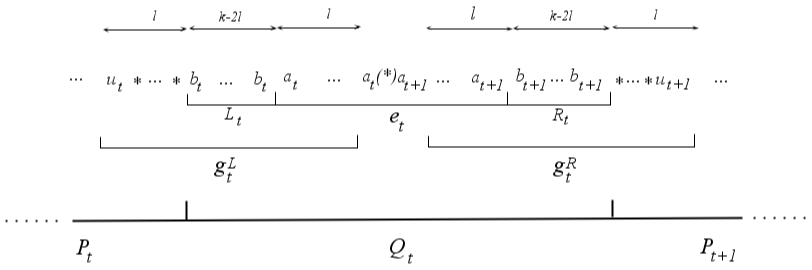}
\end{center}
\caption{Construction of $C$, $t+1\in T_2$.}
\label{cycleC2}
\end{figure}

\begin{figure}
\begin{center}
\includegraphics[scale=0.45]{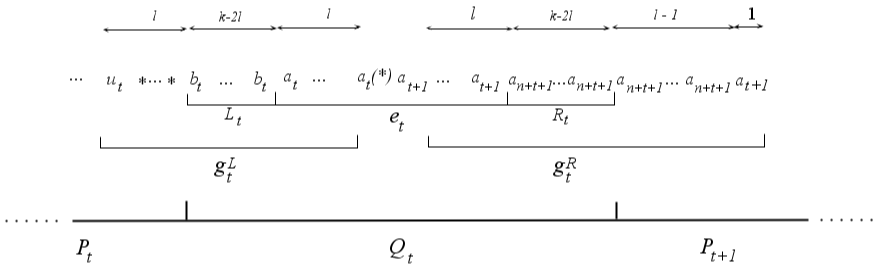}
\end{center}
\caption{Construction of $C$, $t+1\in T_1$.}
\label{cycleC1}
\end{figure}

We will show that the following sequence
\[C= Q_0, P_1, Q_1, P_2, Q_2,  P_3, \dots, Q_{n-1}, P_{n}.\]
spans a hamiltonian $(\ell,k)$-cycle in $H_1\cup H_2+e$. Recall that for each $t\in[n]$, $P_t\subseteq H_2$. Also, each sequence $Q_t$, $t\in[0,n-1]$, consists of a core path ($Q_0'\subseteq H_1\cup H_2+e$ for $t=0$ and just one edge $e_t\in H_1$ for $t\in[n-1]$) and two ``loose ends'' of $k-2\ell$ vertices each. Thus, there are exactly $2n$ edges of $C$ which are not contained in $Q_0\cup P_1\cup\cdots\cup Q_{n-1}\cup P_n$ and require a proof that they also belong to $H_1\cup H_2$. Each of these new edges shares exactly $k-\ell$ vertices with a $Q_t$ and $\ell$ vertices with either $P_t$ ($P_n$ for $t=0$) or $P_{t+1}$, $t=0,\dots,n-1$. Let us denote them by $g_t^L$ and $g_t^R$, respectively (see Figures \ref{cycleC2} and \ref{cycleC1}). For convenience, we set $P_0=P_n$.

Let us first focus on $g_t^L$, $t\in[0,n-1]$. By the construction of $Q_t$ (see \eqref{eee}-\eqref{diagram_q0} for $t=0$ and \eqref{Qt} for $t\ge1$), we have $g_t^L\cap Q_t\subset U_{{t}}$, so $|g_t^L\cap Q_t \cap U_t|=k-\ell$. Further, as  $P_{t}\subset H_2$, among its last $\ell$ vertices there must  be at least one from $U_t$. Since $|g_t^L \cap V(P_{t})|=\ell$, it altogether yields that $g_t^L \in H_2$.
In the same way one can prove that  $g_t^R\in H_2$ for all $t$ such that $t+1\in T_2$ (see Fig. \ref{cycleC2}).

Finally, consider $g_t^R$ with $t+1\in T_1$, (see Fig. \ref{cycleC1}).
  By \eqref{Qt} and (\ref{lpt1}) we have $\{t+1,n+t+1\}\in tr(g_t^R)$, $|g_t^R \cap A_{{t+1}+n}|= k-\ell-1$ and $|g_t^R\cap A_{t+1}|=\ell+1$. Hence,  $g_t^R\in H_1^2$.

\qed
{\color{black}
\section{Concluding remarks}
After fixing an inaccuracy in the first version of our proof, it turned out, quite disappointedly, that Theorem \ref{new}, and thus Corollary \ref{new exact}, does not cover the case $\ell = \lfloor k/2 \rfloor={\color{black} (k-1)/2}$ for odd $k$.
However, a few little changes in the proof can close this gap.
In order to confirm Conjecture \ref{Con} for $\ell = (k-1)/2$, one has to prove Lemmas \ref{pierwszy} and \ref{drugi} for
\[2\leq p = \ell = (k-1)/2, \]
which together will imply a corresponding version of Theorem \ref{new} for $p = \ell = (k-1)/2$, and thus Conjecture \ref{Con} for $\ell=(k-1)/2$.


The change in the proof boils down to replacing $2\lfloor k/2\rfloor +\ell$ with $2\lfloor k/2\rfloor +\ell-1$ in \eqref{mocai} and,  accordingly,
$x_i - 2\lfloor k/2 \rfloor -\ell$ with $x_i - 2\lfloor k/2 \rfloor -\ell+1$ in \eqref{mocbi}. Notice that, for each $i$, $|U_i| = |A_i|+|B_i|$ stays unchanged.
As a result, \eqref{delta} and \eqref{mleqn} remain true, since now $|A_j| \leq 3\lfloor k/2 \rfloor -1$. Moreover, although
inequality \eqref{=l2} is relaxed to
\begin{align*}
\left|\left(V(P_i)\cap A_{j_i}\right)\setminus (e_1 \cup \cdots \cup e_m)\right| \geq \ell,
\end{align*}
it still implies that $deg_{Tr(M)}(j_i) \leq 1$, because $\left|A_{j_i}\right| \leq 2\lfloor k/2 \rfloor +\ell -1$. This  saves
Claims \ref{cl1} and~\ref{cl1b}, while all estimates of the length of $C$ remain intact  (they rely mainly on the
cardinalities of $U_t$ which have not changed). Thus, the proof of Lemma \ref{pierwszy} is retained.

In order to modify the proof of Lemma \ref{drugi}, in Subsection \ref{ij}
one has to choose $i$ and $j$ according to whether $\rho(C_1)\leq \ell - 1$ or $\rho(C_1)\geq \ell$, instead of
 $\rho(C_1)\leq \ell$ or $\rho(C_1)\geq \ell +1$ (and the same for $\rho(C_2)$). This does not affect the structural properties
of the bridge $Q_0$, as consecutive edges intersect in $\ell$ vertices only, but at the same time
strengthens  Fact \ref{ute}  to $|e \cap U_t|\leq \ell-1$. This, in turn,  allows one
to replace the middle part of Fact \ref{pezero} by $|Q_0\cap A_t|\leq \ell-1$, compensating for the decrease of $|A_t|$.

Indeed, since  all  bridges  $Q_1,\dots,Q_{n-1}$, defined in \eqref{Qt},  use together  at most $ 2\lfloor k/2 \rfloor$ vertices from each set $A_t$, $t=2,..., n-1$, this part
of Fact \ref{pezero} implies that there are sufficiently many vertices in the sets $A_t$,  $t=2,..., n-1$, to construct
all bridges (including $Q_0$). On the other hand, for each $t \in \{1,n\}$, the bridges $Q_1,\dots,Q_{n-1}$ require
 only at most $\lfloor k/2 \rfloor$ vertices from $A_t$. Hence, by the first line of Fact \ref{pezero} (with $p=\ell = (k-1)/2$),
we have
\[k-p+\lfloor k/2 \rfloor  = k \leq 3 (k-1)/2  -1      = |A_t|,\]
since $k\geq 5$. Consequently, the construction of all bridges can be completed.
As the remainder  of the proof of Lemma \ref{drugi} does not involve  the (modified) cardinalities of the sets $A_t$,
 the construction of the Hamiltonian cycle $C$ can be finalized basically in the same way as presented in Subsections \ref{conP} and \ref{conC}.

Let us summarize  that, owing to the above extension, Conjecture \ref{Con} is now confirmed for $\ell=1$, all $(k-1)/3\le\ell\le k/2$, and all $\ell\ge0.8k$. We believe that the two missing ranges of $\ell$ will require some new ideas.}

{\color{black}
\section*{Acknowledgements } We would like to thank both referees for several remarks and suggestions leading to a great improvement of the
exposition of the paper. We are especially indebted   to referee X who found an inaccuracy in an earlier version of the proof (c.f. Concluding Remarks).}


\section*{Appendix: Properties of function $\nu$}

In \cite{RZ2} we proved the following simple facts.
\begin{pro}[\cite{RZ2}]\label{obs31} Function $\nu$ has the following properties.
\begin{enumerate}
\item[(a)] For every  $ x\ge (k-3)(k-1)$,
$\quad\nu(x)\geq x+ \left\lfloor \frac{x}{k-1}\right\rfloor+3- k.$
\item[(b)]
For every $x\ge k-2$, $\quad\nu(x)\le kx.$
\item[(c)]
For all $x\ge 2$,$\quad
\nu(x-1)\geq \nu(x)-k.$
\end{enumerate}
\end{pro}
\noindent We will now note three consequences of the above proposition. For $x\geq k^3$ it follows from Proposition \ref{obs31}(a) that
\begin{align}\label{x_nu}
x\leq \frac{k-1}{k}\nu(x)+\frac{(k-1)(k-2)}{k} \leq \frac{k}{k+1}\nu(x)\le \nu(x).
\end{align}
{\color{black}
Indeed, after dropping the floor in (a), we get the first inequality above, while the second inequality is equivalent to $\nu(x)\ge (k+1)(k-1)(k-2)$, which is true by the assumption on $x$.}
Moreover, since $\nu(x)$ equals $\ell$ modulo $k-\ell$,  Proposition \ref{obs31}(c) can be strengthened to yield, for $x\ge2$,
\begin{align}\label{roznica1}
\nu(x)=\nu(x-1)\quad\mbox{or}\quad\nu(x)-\nu(x-1)=k-\ell.
\end{align}
Finally, by iterating the inequality of Proposition \ref{obs31}(c) $t$ times, we have
\begin{align}\label{iterat}
\nu(x+t)\le\nu(x)+tk.
\end{align}


It follows directly from these definitions that
\begin{equation}\label{moreless}
z\geq \nu\left(\mu(z)\right)\quad\mbox{and}\quad z<\nu\left(\mu^*(z)\right).
\end{equation}
The following properties of functions $\nu,\mu$, and $\mu^*$ will turn out to be crucial in our proofs. 

\begin{pro}\label{dwa_mu} We have
\begin{align} \label{roznica}
\nu\left(\mu^*(z)\right)-\nu\left(\mu(z)\right)=k-\ell,
\end{align}
\begin{align}\label{dwa_mu_1}
\nu\left(\mu(z)\right)=\nu\left( \mu(z)-(k-2\ell)\right),\quad\mbox{and}\quad\nu\left(\mu^*(z)\right)=\nu\left( \mu^*(z)+(k-2\ell)\right).
\end{align}

\end{pro}

\proof Equality \eqref{roznica} follows from \eqref{roznica1} and \eqref{mu}.
In order to deduce \eqref{dwa_mu_1}, we first determine an exact formula for function $\nu$ from which it will follow quickly.
Set $\kappa=k-\ell+1$ and $\beta=2k-4\ell+2$ and notice that

\begin{align*}
\max\{\kappa,\beta\}=\begin{cases} \kappa \qquad\text{ if } \ell\geq \frac{k+1}{3}, \\
\beta \qquad\text{ if } \ell < \frac{k+1}{3}.
\end{cases}
\end{align*}
Let us choose an integer $x$ and define  integers $q:=q(x,k,\ell)$ and $r:=r(x,k,\ell)$ by setting
\begin{equation}\label{x}
x-\kappa=q\max\{\kappa,\beta\}+r,\quad\mbox{where}\quad 0\le r\le\max\{\kappa,\beta\}-1.
\end{equation}
We claim that
\begin{equation}\label{nuwzor}
\nu(x)=\begin{cases}
q(2k-2\ell)+k \quad\qquad\qquad\text{ if } r\leq k-2\ell\\
q(2k-2\ell)+2k-\ell \quad\qquad\text{ if } r\geq k-2\ell+1
\end{cases}
\end{equation}

Formula \eqref{nuwzor} shows  that $\nu(x)$ is a step functions which is constant on intervals (steps) of lengths, alternately,  $k-2\ell+1$, and $\max\{\kappa,\beta\}-1-(k-2\ell)\ge \beta-1-(k-2\ell)=k-2\ell+1$. This, together with the definitions of $\mu$ and $\mu^*$, implies equalities \eqref{dwa_mu_1}. Indeed, let, for instance, $x=\mu(z)$ for some $z$. Then $\nu(x)\le z$ but $\nu(x+1)>z$. In view of \eqref{nuwzor} this means that in the expression \eqref{x} we have either $r=k-2\ell$  or $r=\max\{\kappa,\beta\}-1$, that is, $x$ is at the right end of a step of $\nu$. Thus, clearly,  $\nu(x-(k-2\ell))=\nu(x)$, as requiblack. For the second equality in \eqref{dwa_mu_1}, observe that if $x=\mu^*(z)$, then $\nu(x)>z$ but $\nu(x-1)\le z$, so $x$ sits at the left end of a step of $\nu$.


In order to show \eqref{nuwzor}, we will first prove an upper bound valid for all $(\ell,k)$-paths $P$ satisfying (\ref{constr}) and then construct a particular $(\ell,k)$-path $P_0$ which achieves this bound.

Let $P$ be an $(\ell,k)$-path
with  $t$ edges satisfying (\ref{constr}).
Let $e_1,\dots, e_t$ be the edges of $P$ in the linear order underlying $P$.
Set $s=\left\lfloor\tfrac{t+1}2\right\rfloor$.
Clearly, $t\in \{2s-1,2s\}$. Further, set
\begin{align*}
f_i=e_{2i-1}\cup e_{2_i}\setminus e_{2i+1}, \;\;\; i=1,\dots,s-1.
\end{align*}
Since, by (\ref{constr}), $|e_{2i-1}\cap U|\geq \kappa$ for
each $i\in \{1,\dots,s\}$, we have $|f_i\cap U|\geq \kappa$ for
each $i\in \{1,\dots,s-1\}$, too.
However, if $\ell<(k+1)/3$, then this bound can be improved. As, also,  $|e_{2i}\cap U|\geq \kappa$ for
each $i\in \{1,\dots,s-1\}$, we infer that
\begin{align*}
\left|\left(e_{2i}\setminus (e_{2i-1}\cup e_{2i+1})\right) \cap U\right| \geq \kappa-2\ell=k-3\ell+1.
\end{align*}
Therefore,
\begin{align*}
|f_i \cap U|\geq \beta \;\;\; i=1,\dots, s-1,\quad\mbox{and}\quad
|e_{2s-1} \cap U|\geq \kappa.
\end{align*}
  Because $f_1,f_2,\dots,f_{s-1},e_{2s-1}$ are pairwise
disjoint, this implies, in view of \eqref{x},  that $s-1\leq q$.
Also by (\ref{constr}), if $t=2s$, then
$$\left|\left(e_{t}\setminus e_{2s-1}\right)\cap U\right|\geq \kappa-\ell=k-2\ell+1.$$
Thus, if $r\leq k-2\ell$, then  $t=2s-1$ and
\begin{align*}
|V(P)|=\sum_{i=1}^{s-1}|f_i|+|e_{2s-1}|=(s-1)(2k-2\ell)+k \leq q(2k-2\ell)+k.
\end{align*}
Otherwise,   $t\leq 2s$
and
\begin{align*}
&|V(P)|= \sum_{i=1}^{s-1}|f_i|+|e_{2s-1}\cup e_{2s}|=(s-1)(2k-2\ell)+2k-\ell\\
&\leq q(2k-2\ell)+2k-\ell.
\end{align*}

To show equality, let us construct $P_0$ satisfying (\ref{constr}) which achieves this bound.
We will represent $P_0$ as a binary sequence $Q$ over the alphabet $\{u,w\}$, where each vertex
of $U$ is represented by  $u$ and each vertex of $V(P_0)\cap W$ is represented by  $w$
(and the edges of $P_0$ follow the sequence $Q$ according to the definition of an $(\ell,k)$-path).

Assume first that $\ell \geq \frac{k+1}{3}$. Sequence $Q$ consists of $q$ identical blocks plus another block at the end (see diagram (\ref{cdc}) below).
 Each block begins with a $u$-run of length $\kappa-\ell$, followed by a $w$-run of length $\ell-1$, followed by a $u$-run of length $\ell$, followed by a $w$-run of length $k-2\ell$.
 The final block begins with the same runs as all previous blocks, that is, a $u$-run of length $\kappa-\ell$, followed by a $w$-run of length $\ell-1$, followed by a $u$-run of length $\ell$. If $r\le k-2\ell$, then this is it, except that we arbitrarily convert $r$ symbols $w$ to $u$.
If $r\ge k-2\ell+1$, we add a $u$-run of length $r$  followed by a $w$-run of length $k-\ell-r$, creating one more edge. In this case there is no need for any final alteration.
\begin{align}\label{cdc}
\overbrace{\underbrace{u,\dots,u}_{\kappa-\ell}\underbrace{w,\dots,w}_{\ell-1},\underbrace{u,\dots,u}_{\ell}}^{e_1}
\underbrace{w,\dots,w}_{k-2\ell}
\overbrace{\underbrace{u,\dots,u}_{\kappa-\ell}\underbrace{w,\dots,w}_{\ell-1},\underbrace{u,\dots,u}_{\ell}}^{e_3}
\underbrace{w,\dots,w}_{k-2\ell} \nonumber\\
&\cdots \\
\overbrace{\underbrace{u,\dots,u}_{\kappa-\ell}\underbrace{w,\dots,w}_{\ell-1},\underbrace{u,\dots,u}_{\ell}}^{e_{2q-1}}
\underbrace{w,\dots,w}_{k-2\ell}
\overbrace{\underbrace{u,\dots,u}_{\kappa-\ell}\underbrace{w,\dots,w}_{\ell-1},\underbrace{u,\dots,u}_{\ell}}^{e_{2q+1}}
\underbrace{(u,\dots,u}_{r},\underbrace{w,\dots,w)}_{k-\ell-r} \nonumber \\
\nonumber
\end{align}
It is easy to check that $P_0$ satisfies (\ref{constr}). Indeed, the number of symbols $u$ equals $q\kappa+\kappa+r=x$ which agrees with \eqref{x}. Moreover, every edge of $P_0$ covers at least $\kappa$ symbols $u$. This is clearly seen on diagram (\ref{cdc}) for edges  $e_{2i+1}$, $i=1,\dots,q$. However, since $\ell\ge k-\ell$, every edge $e_{2i}$, $i=1,\dots,q$, also contains at least $\kappa-\ell+\ell=\kappa$ symbols $u$. And the last edge, $e_{2q+2}$, if present, contains at least $\ell+r\ge\ell+(k-2\ell+1)=\tau$ symbols $u$ too. (We write ``at least'' as we do not count possible converts from $w$ to $u$.)
Finally, as desiblack (cf. \eqref{nuwzor}),

\begin{equation}\label{VP0}
|V(P_0)|=\begin{cases}
q(k+k-2\ell)+k=q(2k-2\ell)+k \quad\qquad\qquad\qquad\qquad\text{ if } r\leq k-2\ell\\
q(k+k-2\ell)+k+(k-\ell)=q(2k-2\ell)+2k-\ell \quad\qquad\text{ if } r\geq k-2\ell+1.
\end{cases}
\end{equation}

For $\ell<(k+1)/3$ we modify the above construction by replacing each $w$-run of length $k-2\ell$  by a $u$-run of length $k-3\ell+1$ followed by a $w$-run of length $\ell-1$. Again, it is easy to check that both, \eqref{constr} and \eqref{VP0}, hold. Indeed, the total number of symbols $u$ is $q(\kappa+k-3\ell+1)+\kappa+r=q\beta+\kappa+r$ which, again, agrees with \eqref{x}. Moreover, each edge of $P_0$ covers at least $\kappa$ symbols $u$.
Again, this is clear for odd-index edges, while for even indices notice that, this time, $\ell<\kappa-\ell$, so these edges contain each at least $\ell+(k-3\ell+1)+\ell=\tau$ symbols $u$. Finally, the above modification of our construction does not change the total number of vertices in $P_0$, so $|V(P_0)|$ is the same as in \eqref{VP0}.
\qed

\medskip

By \eqref{dwa_mu_1} in Proposition \ref{dwa_mu} and the definitions of $x$ and $x^*$ above,
\begin{equation}\label{mocV_wn}
\nu(x-2\lfloor k/2 \rfloor)=\nu(\mu(z))=\nu(\mu(z)-(k-2\ell))=\nu\left(x-2\lfloor k/2 \rfloor-(k-2\ell)\right)
\end{equation}
and
\begin{equation}\label{mocV_wn1}
\nu(x^*-2\lfloor k/2 \rfloor)=\nu(\mu^*(z)+(k-2\ell))=\nu(\mu^*(z))=\nu\left(x^*-2\lfloor k/2 \rfloor-(k-2\ell)\right).
\end{equation}
Also, by Proposition \ref{obs31}(b), the monotonicity of $\nu$, (\ref{roznica}), (\ref{moreless}), the definition of $z$ in \eqref{mocV}, and (\ref{Nnk}),
\begin{align}\label{lbx}
x&\geq \frac{\nu(x)}{k}\ge\frac{\nu(\mu(z))}k= \frac{\nu(\mu^*(z))-(k-\ell)}{k}>\frac{z-k}k\nonumber \\
&\geq \frac{N}{kn} +\frac{4k^2}n-4\geq 11k^4-4\geq 10k^4.
\end{align}
 In particular, $x-2\lfloor k/2 \rfloor\ge k^3$, which justifies several future applications of \eqref{x_nu}.

On the other hand, by \eqref{x_nu},\eqref{mocV}, \eqref{iterat}, (\ref{moreless}), and (\ref{Nnk}),
\begin{align}\label{ubx}
x&\leq \nu(x)\le\nu(\mu(z)+k)\le\nu(\mu(z))+k^2\le z+k^2 \leq \frac{N+4k^3}n+k^2\le12k^5.
\end{align}

\begin{pro}\label{Pxi}
There exist $x_i \in \{x,x^*\}$, $i=1,\dots,n$, such that
\begin{align}
nz<\sum_{i=1}^n \nu(x_i-2\lfloor k/2 \rfloor)\leq nz+k-\ell. \label{xi}
\end{align}
\end{pro}
\proof Set $y:=\nu\left( x-2\lfloor k/2 \rfloor \right)$ and {\color{black} $y^*=\nu\left( x^*-2\lfloor k/2 \rfloor \right)$}.
By \eqref{mocV}, \eqref{mocV_wn}, and \eqref{mocV_wn1}, $y=\nu(\mu(z))$ and $y^*=\nu(\mu^*(z))$.
 Thus, by (\ref{roznica}), $y^*-y=k-\ell$, and, by \eqref{moreless}, $y\le z$ while $y^*>z$. We are going to show by induction on $m=1,\dots,n$ that there exists a choice of $x_i\in \{x,x^*\}$, $i=1,\dots,m$, such that (\ref{xi}) is satisfied with $n$ replaced by $m$. Indeed, let $x_1=x^*$, then $ z<y^*=y+(k-\ell)\le z+(k-\ell)$. Fix $2\le m\le n$ and assume the statement is true for  $m-1$. Set $\Sigma:=\sum_{i=1}^{m-1} \nu(x_i-2\lfloor k/2 \rfloor)$. Then
 $$mz<\Sigma+y^*\le mz+2(k-\ell),\quad\mbox{while}\quad mz-(k-\ell)<\Sigma+y\le mz+(k-\ell).$$
  Since $(\Sigma+y^*)-(\Sigma+y)=k-\ell$, we have either $\Sigma+y^*\le mz+(k-\ell)$ or $mz<\Sigma+y$, which completes the proof.
\qed

\noindent{\bf Proof of Proposition \ref{PxiI}.} The R-H-S of \eqref{dwie_nier} is the L-H-S of (\ref{xi}).
On the other hand, by the R-H-S of (\ref{xi}), (\ref{x_nu}),  and~(\ref{lbx}),
\begin{align*}
&\sum_{i\in I} \nu(x_i-2\lfloor k/2 \rfloor)\leq\sum_{i=1}^{n} \nu(x_i-2\lfloor k/2 \rfloor)-\min_{1\le i\le n} \nu(x_i-2\lfloor k/2 \rfloor) \\
&\leq N+4k^3-(3k-4\ell)n + k - \nu(x-2\lfloor k/2 \rfloor)\\
&\leq N+4k^3-(3k-4\ell)n+2k-x\leq N-(3k-4\ell)n-8k^4,
\end{align*}
which is the L-H-S of \eqref{dwie_nier}.
\qed

Indeed,  by \eqref{mocai}, \eqref{mocbi}, (\ref{mocV}), (\ref{x_nu}), \eqref{xi}, and (\ref{Nnk})
\begin{align*}
&\sum_{i=1}^{2n}|A_i|+\sum_{i=1}^n|B_i|=\sum_{i=1}^n x_i+n\left( 2k-2\ell-3 \right) = \sum_{i=1}^n (x_i-2\lfloor k/2 \rfloor)
+n\left( 2\lfloor k/2 \rfloor +2k-2\ell-3 \right) \\
&\leq \frac{k}{k+1}\sum_{i=1}^n \nu(x_i-2\lfloor k/2 \rfloor)+3kn <\frac{k}{k+1}\left(N+4k^3-(3k-4\ell)n+2k\right)+3kn \\
&<N-\frac N{k+1}+\frac{k^2}{k+1}\left(4k^2-3n+2\right)+{\color{black} 5kn}<N-\left(\frac N{k+1}-{\color{black}5kn}\right)< N-4k^4n.
\end{align*}
  Thus, for each $i=n+1,\dots,2n$, we have
\begin{equation}\label{lbb}
b_i\ge\left\lfloor\frac1n\sum_{j=n+1}^{2n}b_j\right\rfloor\ge 4k^4.
\end{equation}
{\color{black} while, trivially,}
\begin{align}\label{ubb}
b_i\leq \left\lceil\frac1n\sum_{j=n+1}^{2n}b_j\right\rceil \le N/n+1\leq 12k^5,
\end{align}
where the last inequality follows by \eqref{Nnk}.



\end{document}